\documentclass[10pt,a4paper]{amsart}
\usepackage[english]{babel}
\usepackage[english]{translator}
\usepackage[T1]{fontenc}
\usepackage[utf8]{inputenc}

\usepackage{amsmath,amsthm,amssymb,amsfonts}
\usepackage{bbm}
\usepackage[foot]{amsaddr}

\usepackage{graphicx}
\usepackage{subfigure}
\usepackage{psfrag}

\usepackage{enumitem}
\usepackage{mathtools}

\usepackage[colorlinks=true]{hyperref}
\hypersetup{
  pdffitwindow=false,
  pdfhighlight=/O,
  pdfnewwindow,
  colorlinks=true,
  citecolor=red,            
  linkcolor=blue,             
  menucolor=blue,            
  urlcolor=blue,             
  pdfpagemode=UseOutlines,
  bookmarksnumbered=true,
  linktocpage,
  pdftitle={An IMEX-RK scheme for capturing similarity solutions in the
  multidimensional Burgers' equation},
  pdfsubject={An IMEX-RK scheme for capturing similarity solutions in the
  multidimensional Burgers' equation},
  pdfauthor={Jens Rottmann-Matthes},
  pdfkeywords={similarity solutions, freezing method, relative equilibria,
  stability, imex-rk, central scheme, burgers' equation},
  pdfcreator={pdflatex,bibtex},
  pdfproducer={LaTeX mit hyperref}
}

\usepackage{xcolor}

\usepackage[
  hmarginratio={1:1},     
  vmarginratio={1:1},     
  textwidth=450pt,        
  heightrounded,          
]{geometry}

\newcommand{\1}{{\mathbbm{1}}}
\newcommand{\R}{\ensuremath{\mathbb{R}}}
\newcommand{\Z}{\ensuremath{\mathbb{Z}}}

\newcommand{\bfk}{\ensuremath{\mathbf{k}}}
\newcommand{\bfa}{\ensuremath{\mathbf{a}}}
\newcommand{\K}{\ensuremath{\mathbb{K}}}
\newcommand{\vect}[1]{\begin{pmatrix} #1%
\end{pmatrix}}
\newcommand{\smallvect}[1]{\left(\begin{smallmatrix} #1%
\end{smallmatrix}\right)}
\newcommand{\halb}{{\frac{1}{2}}}
\newcommand{\wt}[1]{\widetilde{#1}}
\newcommand{\wh}[1]{\widehat{#1}}
\DeclareMathOperator{\mm}{mm}

\DeclareMathOperator{\divergence}{div}
\DeclareMathOperator{\gradient}{\nabla}
\DeclareMathOperator{\laplace}{\Delta}
\DeclareMathOperator{\vol}{vol}

\newcommand{\Dxi}{\Delta \xi{}}

\newcommand{\Deta}{\Delta \eta{}}

\newcommand{\Dxihalb}{\Delta \xi_{\frac{1}{2}}}
\newcommand{\Dt}{\Delta t{}}
\newcommand{\Dtau}{\Delta \tau{}}
\newcommand{\tdd}[1]{\tfrac{\partial}{\partial #1}}
\newcommand{\tDD}[1]{\tfrac{d}{d #1}}
\newcommand{\dd}[1]{\frac{\partial}{\partial #1}}
\newcommand{\DD}[1]{\frac{d}{d #1}}

\newcommand{\e}{\mathrm{e}}

\newcommand{\sol}[1]{{\underline{#1}}}
\newcommand{\CC}{\mathcal{C}}
\newcommand{\liealg}[1]{\ensuremath{\mathfrak{#1}}}

\theoremstyle{plain}
\newtheorem*{definition*}{Definition}%
\newtheorem{definition}{Definition}[section]
\newtheorem{theorem}[definition]{Theorem}

\newtheorem*{lemma*}{Lemma}

\newtheorem*{assumption*}{Hypothesis}
\newtheorem{assumption}{Hypothesis}
\theoremstyle{definition}
\newtheorem*{remark*}{Remark}
\newtheorem{remark}[definition]{Remark}
\newtheorem*{remarks*}{Remarks}

\allowdisplaybreaks

\begin{document}
\title[An IMEX-RK scheme for capturing similarity solutions in the
  multidimensional Burgers' equation]{An IMEX-RK scheme for capturing similarity solutions in the multidimensional Burgers' equation}
\setlength{\parindent}{0pt}
\begin{center}
\normalfont\huge\bfseries{\shorttitle}\\
\vspace*{0.25cm}
\end{center}

\vspace*{0.8cm}
\noindent
\hspace*{.33\textwidth}
\begin{minipage}[t]{0.4\textwidth}
\noindent
\textbf{Jens Rottmann-Matthes}\footnotemark[1] \\
Institute for Analysis \\
Karlsruhe Institute of Technology \\
76131 Karlsruhe \\
Germany
\end{minipage}~\\
\footnotetext[1]{e-mail: \textcolor{blue}{jens.rottmann-matthes@kit.edu}, phone: \textcolor{blue}{+49 (0)721 608 41632}, \\
supported by CRC 1173 'Wave Phenomena: Analysis and Numerics', Karlsruhe Institute of Technology}
~\vspace*{0.6cm}

\noindent
\hspace*{.33\textwidth}
Date: \today
\normalparindent=12pt

\vspace{0.4cm}
\begin{center}
\begin{minipage}{0.9\textwidth}
  {\small
  \textbf{Abstract.} 
  In this paper we introduce a new, simple and efficient numerical scheme for the
  implementation of the freezing method for capturing similarity solutions in
  partial differential equations.
  The scheme is based on an IMEX-Runge-Kutta approach for a method of lines
  (semi-)discretization of the freezing partial differential
  algebraic equation (PDAE).  We prove second order convergence for the time
  discretization at smooth solutions in the ODE-sense and we present
  numerical experiments that show second order convergence for the full
  discretization of the PDAE.
  
  As an example serves the multi-dimensional Burgers' equation.  By considering
  very different sizes of viscosity, Burgers' equation can be considered as a
  prototypical example of general coupled hyperbolic-parabolic PDEs.  
  Numerical experiments show that our method works perfectly well for all sizes
  of viscosity, suggesting that the scheme is indeed suitable for capturing
  similarity solutions in general hyperbolic-parabolic PDEs by direct
  forward simulation with the freezing method.
  }
\end{minipage}
\end{center}



\section{Introduction}
Many time-dependent partial differential equations from applications 
exhibit simple patterns.  When these patterns are stable, solutions
with sufficiently close initial data develop these patterns as as time
increases.  They often have important implications on the actual
interpretation of the systems behavior.  A very important and well-known
example of a simple pattern are the traveling wave solutions which
appear in the nerve-axon-equations of Hodgkin and Huxley.  Here they
model the transport of information along the axon of a nerve cell. 

Traveling waves are one of the simplest examples of patterns which are relative
equilibria.  Relative equilibria are solutions to the evolution equation whose
time evolution can be completely described by some curve in a symmetry group
which acts on a fixed profile.  In the case of a traveling wave,
the profile is just the shape of the wave, the curve is a linear function into
the group of the real numbers and the action is the shift of the profile.  Of
course, the slope of the curve is the velocity of the traveling wave.  From an
applications point of view, the velocity of the traveling pulse in the
Hodgkin-Huxley system is of great relevance as it quantifies how fast
information is passed on in this system.  

The above example shows that one is often interested in relative equilibria of
the underlying PDE problem and also in their specific constants of motion like
their velocity.  In a more mathematical way, one has an evolution equation
of the form 
\begin{equation}\label{eq:0.1}
  u_t=F(u),
\end{equation} 
and a traveling wave is a solution $u(x,t)$ of \eqref{eq:0.1} of the form
$u(x,t)=\sol{u}(x-\sol{c}t)$, where $\sol{u}$ is a fixed profile and $\sol{c}$
is its velocity.  When \eqref{eq:0.1} is considered in a co-moving frame, i.e.\
with the new spatial coordinate $\xi=x-\sol{c}t$, the profile $\sol{u}$ becomes
a steady state of the co-moving equation
\begin{equation}\label{eq:0.2}
  v_t=F(v)+\sol{c}v_\xi.
\end{equation}
Similarly, there are evolution equations which exhibit rotating waves or spiral
waves.  In $d=2$ spatial dimensions these are solutions of the form
$u(x,t)=\sol{u}
\bigl(\smallvect{\cos(\sol{c}t)&-\sin(\sol{c}t)\\
\sin(\sol{c}t)&\cos(\sol{c}t)}x\bigr)$.
Considering the evolution equation \eqref{eq:0.1} in a co-rotating frame, i.e.\
with the new spatial coordinates
$\xi=\smallvect{\cos(\sol{c}t)&-\sin(\sol{c}t)\\\sin(\sol{c}t)&\cos(\sol{c}t)}x$,
the profile of a rotating wave becomes a steady state of the co-rotating equation
\begin{equation}\label{eq:0.3}
  v_t=F(v)+\sol{c}\xi_2v_{\xi_1}-\sol{c}\xi_1v_{\xi_2}.
\end{equation}
Again, from an applications point of view, not only the profile $\sol{u}$ is of
interest, but also the velocity $\sol{c}$ with which it rotates.

Typically, the velocities $\sol{c}$ of a traveling or rotating wave (or similar
numbers for other relative equilibria) are not known in advance and thus
the optimal co-moving coordinate frame, in which the solution becomes a
steady state, cannot be used.  A
method, which calculates the solution to the Cauchy problem for \eqref{eq:0.1}
and in parallel a suitable reference frame is
the freezing method, independently introduced in
\cite{BeynThuemmler:2004} and
\cite{RowleyKevrekidisMarsdenLust:2003}, see also
\cite{BeynOttenRottmann:2013}.  A huge advantage of the method is that
asymptotic stability with asymptotic phase of a traveling (or rotating) wave
for the original system becomes asymptotic stability in the sense of Lyapunov of
the waves profile and its velocity for the freezing system, see
\cite{Thuemmler:2008}, \cite{Rottmann:2012a}, \cite{Rottmann:2012c} and
\cite{BeynOttenRottmann:2013}.  As such the method
allows to approximate the profile and its constants of motion by a direct
forward simulation.  The method not only works for
traveling, rotating or meandering waves but also for other relative
equilibria and similarity solutions with a more complicated symmetry group,
e.g.\ \cite{BeynOttenRottmann:2013} or \cite{Rottmann:2016a}. 

The idea of the method is to write the equation in new (time dependent)
coordinates and to split the evolution of the solution into an evolution of the
profile and an evolution in a symmetry group which brings the profile into the
correct position via the group action.  For example, in case of an $E(2)$
equivariance of the evolution equation ($E(2)$ is the Euclidean group of the
plane), when we allow for rotation
and translation in $\R^2$, this leads to 
\begin{subequations}\label{eq:0.freeze}
\begin{equation}\label{eq:0.4}
  v_t=F(v)+\mu_1\partial_{\xi_1}v+\mu_2 \partial_{\xi_2}v+\mu_3
  \bigl(\xi_2\partial_{\xi_1}v-\xi_1\partial_{\xi_2}v\bigr).
\end{equation}
In \eqref{eq:0.4} we have the new time dependent unknowns $v$ and $\mu_1,\mu_2,\mu_3$, where
$\mu_j\in\R$.  To cope with the additional degrees of freedom due to the
$\mu_j$, \eqref{eq:0.4} is supplemented by algebraic equations, so called phase
conditions, which abstractly read
\begin{equation}\label{eq:0.5}
  0=\Psi(v,\mu).
\end{equation}
\end{subequations}
The complete system \eqref{eq:0.freeze} is called the freezing
system and the Cauchy problem for it can be implemented on a computer.  Note
that \eqref{eq:0.freeze} in fact is a partial differential algebraic equation
(PDAE).

In many cases the Cauchy problem for \eqref{eq:0.freeze} can be solved by using
standard software packages like COMSOL Multiphysics, see for example
\cite{BeynOttenRottmann:2013} or \cite{BeynOttenRottmann:2016}.
Nevertheless, in other cases these standard toolboxes may not work at all or may
not give reliable results.  For example this is the case for partly parabolic
reaction-diffusion equations, i.e.\ reaction-diffusion equations in which not
all component diffuse, as is the case in the important Hodgkin-Huxley equations.
In this case the freezing method leads to a parabolic equation that is
nonlinearly coupled to a hyperbolic equation with a time-varying
principal part.  Other examples which cannot easily be
solved using standard packages include hyperbolic conservation laws or coupled
hyperbolic-parabolic PDEs, which appear in many important applied problems.

In this article we present a new, simple and robust numerical discretization
of the freezing PDAE which allows us to do long-term simulations of
time-dependent PDEs and to capture similarity solutions, also for viscous and
even inviscid conservation laws by the freezing method.

We derive our fully discrete scheme in two stages.  First we do a spatial
discretization of the freezing PDAE and obtain a method of lines (MOL) system.
For the spatial discretization we employ a central scheme for hyperbolic
conservation laws from Kurganov and Tadmor \cite{KurganovTadmor:2000}.
Namely, we adapt the semi-discrete scheme derived in
\cite{KurganovTadmor:2000} to the case, when the flux may depend also on the
spatial variable.  This yields a second order semi-discrete central
scheme.  It has the important property that it does not require information
of the local wave structure besides an upper bound on the local wave speed.
In particular, no solutions to Riemann problems are needed.

The resulting MOL system is a huge ordinary differential algebraic
equation (ODAE) system, which has parts with very properties.  On the
one hand, parabolic parts lead for fine spatial discretizations to very
stiff parts in the equation, for which one should employ implicit
time-marching schemes.  On the other hand, a hyperbolic
term leads to a medium stiffness but becomes highly nonlinear due to the spatial
discretization, so an explicit time-marching scheme is preferable.  To
couple these conflicting requirements, we use an implicit-explicit (IMEX)
Runge-Kutta scheme for the time discretization.  Such schemes were
considered in \cite{AscherRuuthSpiteri:1997}, but here we will apply
them to DAE problems.

As an example we consider Burgers' equation,
\begin{equation}
  \label{eq:burgers1d}
  u_t+\left(\tfrac{1}{2}u^2\right)_x=\nu u_{xx},\quad x\in\R,
\end{equation}
which was originally introduce by J.M. Burgers (e.g. \cite{Burgers:1948}) as a
mathematical model of turbulence.
We also consider the multidimensional generalizations of \eqref{eq:burgers1d}
\begin{equation}
  \label{eq:burgersdd}
  u_t+\tfrac{1}{p}\divergence(a|u|^p)=\nu\Delta u,\quad x\in\R^d,
\end{equation}
where $a\in\R^d\setminus\{0\}$ and $p>1$ are fixed.  
Equation~\eqref{eq:burgers1d} and its generalization to $d$ dimensions are among
the simplest truly nonlinear partial differential equations and, moreover, in the
inviscid ($\nu=0$) case they develop shock solutions.  As such they are often
used as test problems for shock capturing schemes, e.g.\ \cite{HuShu:1999}.  But
they are also of interest from a physical point of view as they are special
cases of \emph{the} ``multidimensional Burgers' equation''
\begin{equation}\label{eq:multidburgers}
  \partial_t \vec{u}+(\vec{u}\cdot \gradient)
  \vec{u}=\nu\laplace\vec{u},
\end{equation}
which has applications in different areas of physics, e.g.\ see the review
\cite{BecKhanin:2007}.

There are mainly three reasons for us, why we choose Burgers' equation
as an example.  First of all, it is a very simple (scalar) nonlinear
equation.  But despite its simplicity, it suits as an example for which
the hyperbolic part dominates by choosing $0<\nu <<1$ and it also suits
as an example for which the parabolic part dominates by choosing $\nu >>0$.
And third, as we will see below, the new terms, introduced by the
freezing method, have properties very similar to the terms appearing in the
method of freezing for rotating waves.

The plan of the paper is as follows.
In Section~\ref{sec:2} we derive the continuous freezing method for
Burgers' equation and present the analytic background
In Section~\ref{sec:3} we explain the spatial discretization of the
freezing PDAE and obtain the method of lines ordinary DAE approximation.
In the subsequent step in Section~\ref{sec:4} we then do the
time-discretization of the DAE with our IMEX Runge-Kutta scheme and show
that it is a second order method for the DAE (with respect to the
differential variables).  In the final Section~\ref{sec:5} we present
several numerical results which show that our method and its numerical
discretization are suitable for freezing patterns in equations for which
the parabolic part dominates and also for equations for which the
hyperbolic part dominates.  Because of this we expect our method to be
well suited also for capturing traveling or rotating waves in 
hyperbolic-parabolic coupled problems.


\section{The Continuous Freezing System}
\label{sec:2}
In this section we briefly review the results from \cite{Rottmann:2016a} and
explain the freezing system for Burgers' equation \eqref{eq:burgersdd}.
For the benefit of the reader, though, we present a simplified and direct
derivation without using the abstract language and theory of Lie groups
and Lie algebras.  We refer to \cite{Rottmann:2016a} for the abstract
approach.

\subsection{The Co-Moved System}
To formally derive the freezing system, we assume that the solution $u$ of the
Cauchy problem for \eqref{eq:burgersdd},
\begin{equation}\label{eq:OCauchy}
  \left\{\begin{aligned}
    u_t&=\nu\laplace u-\frac{1}{p}\divergence_x\bigl( a|u|^p\bigr)=:
    F(u),\\
    u(0)&=u_0,
  \end{aligned}
  \right.
\end{equation}
is of the form 
\begin{equation}\label{eq:ansatz}
  u(x,t)=\frac{1}{\alpha\bigl(\tau(t)\bigr)}\,v\left(
  \frac{x-b\bigl(\tau(t)\bigr)}{\alpha\bigl(\tau(t)\bigr)^{p-1}},\tau(t)\right),
\end{equation}
where $v:\R^d\times\R\to\R$, $b=(b_1,\dots,b_d)^\top:\R\to\R^d$, $\alpha:\R\to(0,\infty)$, and
$\tau:[0,\infty)\to[0,\infty)$ are smooth functions and $\dot{\tau}(t)>0$ for
all $t\in[0,\infty)$.
\begin{remarks*}
  \begin{enumerate}[label=\textup{(\roman*)}]
    \item One can interpret the function $\tau$ as a
      transformation of the time $t$ to a new time $\tau$, the action of the
      scalar function $\alpha$ on the function $v$ can be understood as a
      scaling of the function and the space.  Finally, we interpret the meaning
      of $b$ as a spatial shift.
    \item Note that in \cite{Rottmann:2016a} we also allow for a spatial
      rotation.  We actually do not use the rotational symmetry here because it
      does not appear in one and two spatial dimensions.  Moreover, the
      symmetries we consider here are also present in the multi-dimensional
      Burgers' system \eqref{eq:multidburgers}, whereas the rotational
      symmetry is not.
  \end{enumerate}
\end{remarks*}
A simple calculation shows 
\begin{equation}\label{eq:Fsym}
  F(u)(x,t)
  =
  \frac{1}{\alpha^{2p-1}}
  F(v)\left(\xi,\tau\right),
\end{equation}
where $\xi=\tfrac{x-b}{\alpha^{p-1}}$ denotes new spatial coordinates.
Moreover, from \eqref{eq:ansatz}, we obtain with the chain rule
\begin{equation}\label{eq:dt0}
  \begin{aligned}
    \dd{t}u(x,t)&
    =\frac{d}{dt}\left(\tfrac{1}{\alpha\left(\tau(t)\right)}\,v\left(
    \tfrac{x-b\left(\tau(t)\right)}{\alpha\left(\tau(t)\right)^{p-1}},
    \tau(t)\right)\right)\\
    &=-\frac{\alpha'(\tau)}{\alpha^2(\tau)} \dot{\tau}
    v(\xi,\tau)-(p-1)\frac{\alpha'(\tau)}{\alpha^2} \xi^\top \gradient_\xi
    v(\xi,\tau) \dot{\tau}
    -\frac{b'(\tau)^\top}{\alpha(\tau)^p}\gradient_\xi v(\xi,\tau)
    \dot{\tau}+\frac{1}{\alpha(\tau)} v_\tau(\xi,\tau)\dot{\tau}.
  \end{aligned}
\end{equation}
By setting
\begin{equation}\label{eq:rec0}
  \mu_1(\tau)=\frac{\alpha'(\tau)}{\alpha(\tau)}\in\R\quad\text{and}\quad
  \mu_{i+1}(\tau)=\frac{b_i'(\tau)}{\alpha(\tau)^{p-1}}\in\R,i=1,\dots,d,
\end{equation}
\begin{equation}\label{eq:deff}
  \phi_1(\xi,v)=-(p-1)\divergence_{\xi}(\xi v)-(1 -d(p-1))v\quad\text{and}\quad
  \phi_2(\xi,v)
  =-\gradient_\xi v=-\bigl(\dd{\xi_j}v\bigr)_{j=1}^d,
\end{equation}
we can write \eqref{eq:dt0} as
\begin{equation}\label{eq:dt1}
  \dd{t}u(x,t)= \frac{1}{\alpha(\tau)}\Bigl(
  \mu_1(\tau)\phi_1(\xi,v)+
  \mu_{(2:d+1)}(\tau) \phi_2(\xi,v)+v_\tau\Bigr)\dot{\tau},
\end{equation}
where $\mu_{(2:d+1)}=(\mu_2,\dots,\mu_{d+1})$.
Because $u$ solves \eqref{eq:OCauchy}, we obtain from
\eqref{eq:Fsym} and \eqref{eq:dt1} under the assumption that $\tau$
satisfies
\begin{equation}\label{eq:rec0b}
  \dot{\tau}=\alpha(\tau)^{2-2p}.
\end{equation}
in the new $\xi,\tau,v$ coordinates for $v$ the equation
\begin{equation}\label{eq:comove0}
  v_\tau=F(v)-\mu_1 \phi_1(\xi,v)- \mu_{(2:d+1)}\phi_2(\xi,v).
\end{equation}
The key observation for the freezing method is Theorem~4.2 from
\cite{Rottmann:2016a}, which
relates the solution of the original Cauchy problem \eqref{eq:OCauchy}
to the solution of the Cauchy problem for \eqref{eq:comove0} in the new
$(\xi,\tau)$ coordinates.  Using the spaces $X=L^2(\R^d)$ and $Y_1=\bigl\{v\in
H^2(\R^d):\xi^\top\gradient_\xi v \in L^2(\R^d)\bigr\}$ the result from
\cite{Rottmann:2016a} can be rephrased as follows:
\begin{theorem}
  A function $u\in\CC([0,T);Y_1)\cap\CC^1([0,T);X)$ solves \eqref{eq:OCauchy}
  if and only if
  the functions $v\in\CC([0,\wh{T});Y_1)\cap \CC^1([0,\wh{T});X)$,
  $\mu_i\in\CC([0,\wh{T});\R)$, $i=1,\dots,d+1$,
  $\alpha\in\CC^1([0,\wh{T});(0,\infty))$, $b\in \CC^1([0,\wh{T});\R^d)$, and
  $\tau\in \CC^1([0,T);[0,\wh{T}))$ solve the system
  \begin{equation}\label{eq:CoCauchyO}
    \begin{aligned}
      v_\tau&=F(v)-\mu_1 \phi_1(\xi,v)-
      \mu_{(2:d+1)} \phi_2(\xi,v),&v(0)&=u_0,\\
      \alpha_\tau&=\mu_1\,\alpha,&\alpha(0)&=1,\\
      b_\tau&=\alpha^{p-1}\mu_{(2:d+1)}^\top,&b(0)&=0,\\
      \dd{t}\tau&=\alpha(\tau)^{2-2p},&\tau(0)&=0,
    \end{aligned}
  \end{equation}
  and $u$ and $v,\alpha,b,\tau$ are related by \eqref{eq:ansatz}.
\end{theorem}

\subsection{The Freezing System}
To cope with the $d+1$ additional degrees of freedom, due to
$\mu_1,\dots,\mu_{d+1}$, we complement \eqref{eq:CoCauchyO} with
$d+1$ algebraic equations, so called phase conditions, see
\cite{BeynThuemmler:2004}.  

\textbf{Type 1: Orthogonal phase condition.}
We require that $v$, and $\mu_1,\dots,\mu_{d+1}$ which solve
\eqref{eq:comove0}, are chosen such that at each time instance
$\|v_\tau\|_{L^2}^2$ is minimized with respect to $\mu_1,\dots,\mu_{d+1}$.
Therefore, we have
\[ 0=\frac{1}{2}\frac{d}{d\mu} \bigl\| F(v)-\mu_1 \phi_1(\xi,v)-
\mu_{(2:d+1)} \phi_2(\xi,v)\bigr\|_{L^2}^2\]
which is equivalent to 
\begin{equation}\label{eq:ophase}
  \left\{
  \begin{aligned}
    0&=\Bigl\langle \phi_1(\xi,v),F(v)-\mu_1 \phi_1(\xi,v)-
    \mu_{(2:d+1)} \phi_2(\xi,v)\Bigr\rangle_{L^2},\\
    0&=\Bigl\langle \dd{\xi_j}v,F(v)-\mu_1 \phi_1(\xi,v)-
    \mu_{(2:d+1)} \phi_2(\xi,v)\Bigr\rangle_{L^2},\;j=1,\dots,d,
  \end{aligned}
  \right.
\end{equation}
where $\phi_1$ and $\phi_2$ are given in \eqref{eq:deff}.
We abbreviate \eqref{eq:ophase} as $0=\Psi^\text{orth}(v,\mu)$.

\textbf{Type 2: Fixed phase condition.}
The idea of the fixed phase condition is, to require that the
$v$-component of the solution always lies in a fixed,
$d+1$-co-dimensional hyperplane, which is given as the
level set of a fixed, linear mapping.
Here we assume that a ``suitable'' reference function $\wh{u}$ is given and
the $v$-component of the solution always satisfies
\begin{equation}\label{eq:fphase}
  \left\{
  \begin{aligned}
    0&=\Bigl\langle \phi_1(\xi,\wh{u}),\wh{u}-v\Bigr\rangle_{L^2},\\
    0&=\Bigl\langle\dd{\xi_j}\wh{u},\wh{u}-v\Bigr\rangle_{L^2},\;j=1,\dots,d,
  \end{aligned}
  \right.
\end{equation}
where $\phi_1$ and $\phi_2$ are given in \eqref{eq:deff}.
We abbreviate \eqref{eq:fphase} as $0=\Psi^\text{fix}(v)$.

We augment system
\eqref{eq:CoCauchyO} with one of the phase conditions
\eqref{eq:ophase} or \eqref{eq:fphase} and obtain
\begin{subequations}\label{eq:CoCauchy}
  \begin{align}
    \label{eq:CoCauchyA1}
    v_\tau&=F(v)-\mu_1 \phi_1(\xi,v)-\mu_{(2:d+1)} \phi_2(\xi,v),&v(0)&=u_0,\\
    \label{eq:CoCauchyA2}
    0&=\Psi(v,\mu),&&\\
    \label{eq:CoCauchyB}
    \alpha_\tau&=\mu_1\alpha,\;b_\tau=\alpha^{p-1}
    \mu_{(2:d+1)}^\top,&\alpha(0)&=1\in\R,\;b(0)=0\in\R^d,\\
    \label{eq:CoCauchyC}
    \frac{d}{d\tau} t&=\alpha(\tau)^{2p-2},
    &t(0)&=0,
  \end{align}
\end{subequations}
where $\Psi(v,\mu)$ is either $\Psi^{\mathrm{orth}}$ from
\eqref{eq:ophase} or $\Psi^{\mathrm{fix}}$ from \eqref{eq:fphase}.

\begin{remark}\label{rem:3.2}
  \begin{enumerate}[label=\textup{(\roman*)}]
    \item
     Observe that \eqref{eq:CoCauchy} consists of the PDE
     \eqref{eq:CoCauchyA1},
     coupled to a system of ordinary differential equations
     \eqref{eq:CoCauchyB} and \eqref{eq:CoCauchyC}, and coupled to a
     system of algebraic equations \eqref{eq:CoCauchyA2}.  Moreover, in
     \eqref{eq:CoCauchyA1} the hyperbolic part dominates for $0<\nu<<1$ and
     the parabolic part dominates for $\nu>>0$.  Finally,
     note that the ordinary differential equations \eqref{eq:CoCauchyB}
     and \eqref{eq:CoCauchyC} decouple from
     \eqref{eq:CoCauchyA1} and \eqref{eq:CoCauchyA2} can be
     solved in a post-processing step.
   \item Note that the $\phi_1$-term in \eqref{eq:deff} resembles the
     generator of rotation, cf. \eqref{eq:0.4}.
   \item
     Under suitable assumptions on the solution and the reference function, 
     the system \eqref{eq:CoCauchyA1}, \eqref{eq:CoCauchyA2} is a PDAE
     of ``time-index'' 1 in the case $\Psi=\Psi^\text{orth}$ and of
     ``time-index'' 2 in the case $\Psi=\Psi^\text{fix}$, where the index is
     understood as differentiation index, see \cite{MartinsonBarton:2000}.
 \end{enumerate}
\end{remark}


\section{Implementation of the Numerical Freezing Method I: Spatial
Discretization}\label{sec:3}
In this section we derive a spatial semi-discretization (Method of Lines system)
of the freezing partial differential algebraic evolution equation system
\eqref{eq:CoCauchy}.  First we separately consider the
PDE-part \eqref{eq:CoCauchyA1} of the equation on the full domain.  Afterwards
we consider the case of a bounded domain with artificial no-flux boundary
conditions.  Finally, we consider the spatial discretization of the
(low-dimensional) remaining equations
\eqref{eq:CoCauchyA2}--\eqref{eq:CoCauchyC}.
\subsection{Spatial Semi-Discretization of the PDE-Part}\label{sec:3.1}
Recall from \eqref{eq:deff}, that $\mu_1 \phi_1(\xi,v)$ is of
the form
\[\mu_1 \phi_1(\xi,v)=\mu_1 \sum_{j=1}^d\DD{\xi_j} f_{1,j}(\xi,v)-\mu_1
f_{1,0}(v),\]
where $f_{1,0}(v)=\bigl(1-d(p-1)\bigr) v$ and $f_{1,j}(\xi,v)=-(p-1)\xi_j v$.
Note that
$\frac{d}{d\xi_j}f(\xi,v)=\dd{\xi_j}f+\dd{v}f\dd{\xi_j}v$.
Similarly, the term $\mu_{(2:d+1)} \phi_2(\xi,v)$ can be written as
\[\mu_{(2:d+1)} \phi_2(\xi,v)= -\sum_{i=1}^d \mu_{1+i} \dd{\xi_i} v = \sum_{i=1}^d
\mu_{1+i}\sum_{j=1}^d \DD{\xi_j} f_{1+i,j}(\xi,v),\]
where $f_{1+i,j}(\xi,v)=-v$ if $i=j$ and $f_{1+i,j}(\xi,v)=0$ if $i\ne
0$, $i,j=1,\dots,d$.
Furthermore, by also setting $f_{0,j}(v)=\tfrac{1}{p}a_j |v|^p$,
equation~\eqref{eq:CoCauchyA1} can be recast as
\begin{equation}\label{eq:MOL.01}
  v_\tau+\sum_{i=1}^{d+1} \mu_i\Bigl[ \sum_{j=1}^d
  \frac{d}{d \xi_j} f_{i,j}(\xi,v)\Bigr]+\sum_{j=1}^d
  \frac{d}{d\xi_j}f_{0,j}(v)=\sum_{j=1}^d
  \frac{d}{d\xi_j}Q_j(\tdd{\xi_j}v)+\mu_1 f_{1,0}(v).
\end{equation}

For the spatial semi-discretization we adapt the 2nd order semi-discrete
central scheme for hyperbolic conservation laws from
\cite{KurganovTadmor:2000} to our situation.  The details of the
adaptation to space dependent hyperbolic conservation laws of the form 
\begin{equation}\label{eq:MOL.03}
  v_\tau+\sum_{j=1}^d \frac{d}{d \xi_j} f_j(\xi,v)=0
\end{equation}
is presented in Appendix~\ref{sec:KT}.  As is noted in
\cite[\S~4]{KurganovTadmor:2000}, diffusive flux terms and zero order
terms, which appear in \eqref{eq:MOL.01}, can easily
be appended to the semi-discretization of \eqref{eq:MOL.03} by simple
second order finite difference approximations.

For the actual discretization we choose a uniform spatial grid in each
coordinate direction
\[\xi^j_{k-\halb}=(k-\halb)\Delta \xi_j,\;k\in\Z,\;j=1,\ldots,d,\] 
and obtain the rectangular cells (finite volumes)
$C_\bfk:=C_{k_1,\ldots,k_d}:=
\bigtimes_{j=1}^d(\xi^j_{k_j-\halb},\xi^j_{k_j+\halb})\subset\R^d$,
with centers
$\xi_\bfk:=(\xi^1_{k_1},\ldots,\xi^d_{k_d})\in\R^d$, where
$\bfk=(k_1,\ldots,k_d)\in\K$ with index set $\K=\Z^d$.  In the following we interpret
for each
$\bfk\in\K$, $v_\bfk(\tau)$ as an approximation of the
cell-average of $v$ in the cell $C_\bfk$ at time $\tau$,
\[v_\bfk(\tau)\approx\frac{1}{\vol(C_\bfk)}\int_{C_\bfk}v(\xi,\tau)\,d\xi,\quad
  \vol(C_\bfk)=\prod_{j=1}^d \Delta \xi_j.\]
Then the method-of-lines system for \eqref{eq:MOL.01} is given by 
\begin{multline}\label{eq:MOL.04}
    v_\bfk'=-\sum_{j=1}^d \frac{H_{\bfk^j+\halb}^{0,j}
      -H_{\bfk^j-\halb}^{0,j}}{\Delta \xi_j}
    -\sum_{i=1}^{d+1}\mu_i\sum_{j=1}^d\frac{
      H^{i,j}_{\bfk^j+\halb}-H^{i,j}_{\bfk^j-\halb}}{ \Delta \xi_j}
    +\mu_1 f_{1,0}(v_\bfk)\\
    +\sum_{j=1}^d \frac{P^j_{\bfk^j+\halb}-P^j_{\bfk^j-\halb}}{
    \Delta \xi_j} =: F_\bfk(v_{\Z^d},\mu), \quad \bfk\in \K.
\end{multline}
In \eqref{eq:MOL.04} we use the following
notations and abbreviations:\\
For $\bfk=(k_1,\ldots,k_d)\in\K$, we denote 
$\bfk^j\pm\halb:=\bigl(k_1,\dots,k_{j-1},k_j\pm\halb,k_{j+1},\dots,k_d\bigr)$.
Then $H_{\bfk^j\pm\halb}^{i,j}(t)$ is an approximation of the
\emph{``hyperbolic flux''} through the boundary
face
\[(\xi_{k_1-\halb}^1,\xi_{k_1+\halb}^1)\times\dots\times
\{\xi^j_{k_j\pm\halb}\}\times
\dots\times(\xi_{k_d-\halb}^d,\xi_{k_d+\halb}^d)=:\partial
C_{\bfk^j\pm\halb}\]
of $C_\bfk$ due to the flux function $f_{i,j}$, $i=0,\dots,d+1$,
$j=1,\dots,d$.  Moreover, 
we add the regularizing part of the KT-scheme completely to the
$H^{0,j}$-terms.  Namely the term
$\tfrac{a}{2\Delta\xi}\bigl(u_{k-\halb}^- - u_{k-\halb}^+ - u_{k+\halb}^-
+u_{k+\halb}^+\bigr)$ from \eqref{eq:semidiscKT} is added to the discretization
of $\tDD{\xi}f_{0,1}(v)$ in the 1-dimensional case.  In the 2-dimensional case
the term
$\tfrac{a}{2\Delta\xi}\bigl(u_{k_1-\halb,k_2}^- - u_{k_1-\halb,k_2}^+ -
u_{k_1+\halb,k_2}^- +u_{k_1+\halb,k_2}^+\bigr)$ from \eqref{eq:semidiscKT2d} is
added to the discretization of $\tDD{\xi_1}f_{0,1}(v)$ and the term
$\tfrac{a}{2\Delta\eta}\bigl(u_{k_1,k_2-\halb}^- - u_{k_1,k_2-\halb}^+ -
u_{k_1,k_2+\halb}^- +u_{k_1,k_2+\halb}^+\bigr)$ from \eqref{eq:semidiscKT2d} is
added to the discretization of $\tDD{\xi_2}f_{0,2}(v)$.
Hence the $H^{i,j}$, $i=0,\dots,d+1$, $j=1,\dots,d$ take the form
\begin{equation}\label{eq:MOL.05}
  H_{\bfk^j\pm\halb}^{i,j}(t)=\frac{f_{i,j}( \xi_{\bfk^j\pm\halb},
    v_{\bfk^j\pm\halb}^+)+ f_{i,j}( \xi_{\bfk^j\pm\halb},
    v_{\bfk^j\pm\halb}^-)}{2}-
  \delta_{i,0}\bfa\bigl(v_{\bfk^j\pm\halb}^+-v_{\bfk^j\pm\halb}^-\bigr).
\end{equation}
In \eqref{eq:MOL.05} the point of evaluation is
\begin{equation}\label{eq:MOL.05a}
  \xi_{\bfk^j\pm\halb}=\bigl(\xi_{k_1},\ldots,\xi_{k_d}\bigr)\pm\tfrac{1}{2}
  \e_j\Delta \xi_j\in\partial C_\bfk,
\end{equation}
and $v_{\bfk^j+\halb}^+$,
$v_{\bfk^j+\halb}^-$ denote the limits ``from above'' and ``from
below'' at the boundary face $\partial C_{\bfk^j+\halb}$ of the
piecewise linear minmod-reconstruction $\hat{v}$ of
$(v_\bfk)_{\bfk\in\Z^d}$, explained in the
Appendix~\ref{sec:KT}, i.e.
\begin{equation}\label{eq:MOL.05b}
  v_{\bfk^j+\halb}^+=\lim_{h\searrow \halb} \hat{v}(x_\bfk+h\Delta
  x_j\e_j),\; 
  v_{\bfk^j+\halb}^-=\lim_{h\nearrow \halb} \hat{v}(x_\bfk+h\Delta
  x_j\e_j).
\end{equation}
For explicit formulas of $v_{\bfk^j+\halb}^\pm$ in the 1d and 2d case
see \eqref{eq:minmodlim1d} and \eqref{eq:minmodlim2d}, respectively.

Moreover, the value $\bfa$, appearing in $H^{0,j}_{\bfk^j\pm\halb}$,
denotes an upper bound for the maximal wave speeds of the hyperbolic
terms.  The choice of the value $\bfa$ is only formal at the moment,
because such a number does not exist in case of a non-trivial scaling.
But note that for the actual
numerical computations performed in Section~\ref{sec:5}, we have to restrict to
a bounded domain and on bounded domains the number $\bfa$ is always finite.

The term $f_{1,0}(v_\bfk)$ in \eqref{eq:MOL.04} is an approximation of the
average source term $f_{1,0}(v)$ on the cell $C_\bfk$.

Finally, $P_{\bfk^j+\halb}^j$, $\bfk=(k_1,\ldots,k_d)\in\K$, denotes an
approximation of the ``\emph{diffusion flux}'' through the
boundary face $\partial C_{\bfk^j+\halb}$, due to the flux function
$Q_j(\tdd{\xi_j}v)$.
This is simply approximated by
\begin{equation}\label{eq:MOL.05c}
  P_{\bfk^j+\halb}^j=Q_j\left(\tfrac{v_{\bfk+\e_j}-v_\bfk}{\Delta
      \xi_j}\right).
\end{equation}
There is no difficulty in generalizing to diffusion fluxes
of the form $Q_j(v,\tdd{\xi_j}v)$, see
\cite[\S~4]{KurganovTadmor:2000}.
\begin{remark*}
  \begin{enumerate}[label=\textup{(\roman*)}]
  \item The upper index $j$ in \eqref{eq:MOL.04}--\eqref{eq:MOL.05c},
    corresponds to the
    $j$th coordinate direction, i.e.
    \[-\frac{H_{\bfk^j+\halb}^{0,j}
      -H_{\bfk^j-\halb}^{0,j}}{\Delta \xi_j}\\
    -\sum_{i=1}^{d+1}\mu_i\frac{
      H^{i,j}_{\bfk^j+\halb}-H^{i,j}_{\bfk^j-\halb}}{ \Delta \xi_j} +
    \frac{P^j_{\bfk^j+\halb}-P^j_{\bfk^j-\halb}}{ \Delta \xi_j}\] is
    the numerical flux in the $j$th direction and $H$ stands for the
    hyperbolic flux and $P$ for the dissipative flux.
  \item In the derivation of \eqref{eq:MOL.04} we make essential use
    of the the property that the method \eqref{eq:semidiscKT} and
    its multi-dimensional variants (e.g.\ \eqref{eq:semidiscKT2d})
    depend linearly on the flux
    function, see Remark~\ref{rem:KT}.
  \end{enumerate}
\end{remark*}

\subsection{Artificial No-Flux Boundary Conditions}\label{sec:3.2}
Equation \eqref{eq:MOL.04} is an infinite dimensional system of
ordinary differential equations and hence not implementable on a
computer.  Therefore, we restrict \eqref{eq:CoCauchy} to a bounded
domain.  For simplicity we only consider rectangular domains of the
form
\begin{equation}\label{eq:MOL.02}
  \Omega = \bigl\{\xi=(\xi_1,\ldots,\xi_d)\in\R^d: R_j^-\leq \xi_j\leq R_j^+,
  j=1,\ldots,d\bigr\}
\end{equation}
with $R_j^-<R_j^+\in\R$.  For the PDE-part \eqref{eq:CoCauchyA1} we
then impose no-flux boundary conditions on $\partial \Omega$.  We perform
the method-of-lines approach as presented in the previous
Subsection~\ref{sec:3.1}.  For this we assume that the
grid in the $j$th coordinate direction, $j=1,\dots,d$, is given by
\[\xi^j_{k_j-\halb}=R_j^-+k_j\Delta \xi_j,\;k_j=0,\ldots,N_j+1,\;
\text{with}\;\xi^j_{N_j+\halb}=R_j^+.\]
Again $\K:=\{\bfk=(k_1,\dots,k_d)\in\Z^d:
0\leq k_j\leq N_j+1,\,j=1,\dots,d\}$ denotes the index set and 
the cells (finite volumes) are
\[C_\bfk:=C_{k_1,\ldots,k_d}:=
  \times_{j=1}^d(\xi_{k_j-\halb}^j,\xi_{k_j+\halb}^j),\quad\bfk\in\K.\]
The cell $C_\bfk$ has the center $\xi_\bfk=(\xi_{k_1}^1,\ldots,\xi_{k_d}^d)$.

For the resulting discretization of \eqref{eq:CoCauchyA1} on $\Omega$, which is
of the form \eqref{eq:MOL.04}, it is easy to implement no-flux boundary
conditions by imposing for $i=0,\dots,d+1$, $j=1,\dots,d$
\begin{equation}\label{eq:MOL.05d}
  H^{i,j}_{\bfk^j\pm\halb}=0\;\text{ and }\;  P^j_{\bfk^j\pm\halb}=0\;
  \text{ if }\;\xi_{\bfk^j\pm\halb}\in\partial \Omega.
\end{equation}
Moreover, an explicit upper bound $\bfa$ of the
local wave speeds in \eqref{eq:MOL.05} is
\begin{equation}\label{eq:MOL.05e}
  \bfa= \max_{\xi\in \Omega}\max_j
  \Bigl[ \sum_{i=1}^{\dim\liealg{g}}
  \mu_i \dd{v} f_{i,j}(\xi,v)+\dd{v}f_{0,j}(v)\Bigr].
\end{equation}

Thus we obtain a method-of-lines approximation of \eqref{eq:CoCauchyA1} on
$\Omega$ subject to (artificial) no-flux boundary conditions on
$\partial\Omega$
by the formula \eqref{eq:MOL.04} for $\bfk\in\K$
with \eqref{eq:MOL.05}, \eqref{eq:MOL.05a}, \eqref{eq:MOL.05b}, \eqref{eq:MOL.05c},
\eqref{eq:MOL.05d}, and \eqref{eq:MOL.05e}.

\subsection{Spatial Semi-Discretization of ODE- and Algebraic-Part}\label{sec:3.3}
Because \eqref{eq:CoCauchyB} and \eqref{eq:CoCauchyC} are already
(low-dimensional) ordinary differential equations and do not explicitly
depend on the values of the function $v$, nothing has to be done
for their spatial discretizations.  Thus it remains to discretize 
\eqref{eq:CoCauchyA2}, the phase conditions.  In this article we assume
that they are given by one of the formulas  \eqref{eq:ophase} or
\eqref{eq:fphase}.  These integral conditions can easily be discretized
by using average values of the function on a cell.
Namely, we approximate \eqref{eq:ophase} by the system
\begin{equation}\label{eq:MOL.AlgO}
  0=  \sum_{\bfk\in\K} \vol(C_\bfk)
  \Bigl[\left(-\sum_{j=1}^d
    \frac{H^{i,j}_{\bfk^j+\halb}-H^{i,j}_{\bfk^j-\halb}}{\Delta \xi_j}
    +\delta_{i,1}f_{1,0}(v_\bfk)\right)\cdot F_\bfk(v_\K,\mu)\Bigr],\;
  i=1,\ldots,d+1
\end{equation}
of $d+1$ algebraic equations.  Similarly, we
approximate \eqref{eq:fphase} by
\begin{equation}\label{eq:MOL.AlgF}
  0= \sum_{\bfk\in\K} \vol(C_\bfk)
  \Bigl[\left(-\sum_{j=1}^d
    \frac{H^{i,j}_{\bfk^j+\halb}(\wh{u}_\K)-H^{i,j}_{\bfk^j-\halb}
      (\wh{u}_\K)}{\Delta \xi_j} 
    +\delta_{i,1}f_{1,0}(\wh{u}_\bfk)\right)\cdot
  \bigl(v_\bfk-\wh{u}_\bfk\bigr)\Bigr],\; 
  i=1,\ldots,d+1,
\end{equation}
where $\wh{u}_\bfk$ is an approximation of the average value of the reference
function $\wh{u}$ in the cell $C_\bfk$, $\wh{u}_\K=(\wh{u}_\bfk)_{\bfk\in\K}$,
and $H^{i,j}_{\bfk^j\pm\halb}(\wh{u}_\K)$ is given by \eqref{eq:MOL.05} with $v$
replaced by $(\wh{u}_\bfk)_{\bfk\in\K}$.
\begin{remark}
  The evaluation of the phase conditions \eqref{eq:MOL.AlgO} and
  \eqref{eq:MOL.AlgF} is cheap and easily obtained.  
  Because the terms $H^{i,j}_{\bfk^j\pm\halb}$ and
  $f_{1,0}(v_\bfk)$ in \eqref{eq:MOL.AlgO} are anyway needed for
  the evaluation of $F_\bfk$. 
  
  Similarly, for the discretization of
  \eqref{eq:MOL.AlgF} one only needs to calculate the
  $H^{i,j}_{\bfk^j\pm\halb}(\wh{u})$ once and has to remember the result of this
  calculation.  In the special case when one chooses some previous value $v_\K$
  as reference solution, the value of $H^{i,j}_{\bfk^j\pm\halb}(\wh{u})$
  is even already known. 
\end{remark}

\subsection{The Final Method of Lines System}
Recollecting the above discussion, the method of lines approximation 
of the freezing PDAE \eqref{eq:CoCauchy} takes the final form
\begin{subequations}\label{eq:HDAE}
  \begin{align}
    \label{eq:HDAE1}
    v_\bfk'&=-H^0_\bfk(v_\K)-H^1_\bfk(v_\K)\mu
    +P_\bfk(v_\K)=:F_\bfk(v_\K,\mu),\quad
    \bfk\in\K,\\
    \label{eq:HDAE2}
    0&=\begin{cases} 
      H^1(v_\K)^\top
      F(v_\K,\mu),\quad&\text{orthogonal Phase},\\
      H^1(\wh{u}_\K)^\top 
      (v_\K-\wh{u}_\K),\quad&\text{fixed Phase},
    \end{cases}\\
    g'&=r_\text{alg}(g,\mu),\label{eq:HDAE3}\\
    t'&=r_\text{time}(g),\label{eq:HDAE4}
  \end{align}
\end{subequations}
with initial data
$v_\bfk\approx\frac{1}{\vol(C_\bfk)}\int_{C_\bfk}u_0(\xi)\,d\xi$ for
$\bfk\in\K$, $g(0)=\bigl(\alpha(0),b(0)^\top\bigr)^\top=(1,0,\dots,0)^\top$, 
and $t(0)=0$.  In \eqref{eq:HDAE} 
we use the abbreviations
$\mu=(\mu_1,\dots,\mu_{d+1})^\top\in\R^{d+1}$ and
\begin{equation*}
  \begin{aligned}
    H^0_\bfk(v_\K)&=\sum_{j=1}^d \frac{H_{\bfk^j+\halb}^{0,j}
    -H_{\bfk^j-\halb}^{0,j}}{\Delta \xi_j},&
    H^0(v_\K)&=\bigl(H^0_\bfk(v_\K)\bigr)_{\bfk\in\K},\\
    H^1_\bfk(v_\K)&=\left(-f_{1,0}(v_\bfk)+\sum_{j=1}^d
    \frac{ H^{1,j}_{\bfk^j+\halb}-H^{1,j}_{\bfk^j-\halb}}{ \Delta \xi_j},\dots,
    \sum_{j=1}^d\frac{ H^{d+1,j}_{\bfk^j+\halb}-
    H^{d+1,j}_{\bfk^j-\halb}}{ \Delta \xi_j}
    \right),&H^1(v_\K)&=\bigl(H^1_\bfk(v_\K)\bigr)_{\bfk\in\K},\\
    P_\bfk(v_\K)&=
    \sum_{j=1}^d \frac{P^j_{\bfk^j+\halb}-P^j_{\bfk^j-\halb}}{ \Delta
    \xi_j},& P(v_\K)&=\bigl(P_\bfk(v_\K)\bigr)_{\bfk\in\K},\\
    F(v_\K,\mu)&=\bigl(F_\bfk(v_\K,\mu)\bigr)_{\bfk\in\K},\\ 
    g&=\vect{\alpha\\b},\quad
    r_\text{alg}(g,\mu)=\vect{g_1\mu_1\\g_1^{p-1}\mu_{(2:d+1)}^\top}=\vect{\alpha
      \mu_1\\ \alpha^{p-1}\mu_{(2:d+1)}^\top},\\
    r_\text{time}(g)&= g_1^{2p-2}=\alpha^{2p-2},
  \end{aligned}
\end{equation*}
with all terms explained in Section~\ref{sec:3.1}.  Note, that we included
the $f_{1,0}$-term in the first column of $H^1_\bfk$.
We also impose \eqref{eq:MOL.05d} to implement no-flux boundary conditions.
As before $v_\K=(v_\bfk)_{\bfk\in\K}$ is replaced by
$\wh{u}_\K=(\wh{u}_\bfk)_{\bfk\in\K}$ in $H^{i,j}_{\bfk^j\pm\halb}$ and
$f_{1,0}$ for the
definition of $H^1(\wh{u}_\K)$.  
Equation \eqref{eq:HDAE} is written in matrix times vector from, where
we understand $H^0(v_\K)$, $f^1(v_\K)$, and $P(v_\K)$ as vectors in $\R^{\#(\K)}$ and
$H^1(v_\K)$ as a matrix in $\R^{\#(\K),d+1}$, with $\#(\K)$ being the number of
elements in the index set $\K$.


\section{Implementation of the Numerical Freezing Method II: Time
Discretization}\label{sec:4}

In this Section we introduce a new time-discretization of the DAE
\eqref{eq:HDAE}.  First we motivate and explain our scheme and in
then we analyze the local truncation error introduced by it.
Note that for a full convergence analysis the
PDE-approximation error has to be analyzed and one has to cope
with the unbounded operators.  We finish this section by a description
of how the implicit equations can efficiently be solved.

\subsection{An Implicit-Explicit Runge-Kutta Scheme}
\label{sec:imex}
Originating from the structure of the original
problem \eqref{eq:CoCauchyA1}, see
Remark~\ref{rem:3.2}, the ODE-part \eqref{eq:HDAE1} of \eqref{eq:HDAE}
has very different components:

The $P(v_\K)$-part is linear in $v_\K$, but as
the discretization of a (negative) elliptic operator has spectrum which
extends far into the left complex half plane.  Therefore, the
``sub-problem'' $v_\K'=P(v_\K)$ becomes very stiff for fine spatial
grids and an efficient numerical scheme
requires an implicit time discretization.

On the other hand, the terms $H^0(v_\K)$ and $H^1(v_\K)\mu$ are highly
nonlinear in the argument $v_\K$ due to the nonlinearity $f$ and the
reconstruction of a piecewise linear function from cell averages,
see \eqref{eq:MOL.05}, \eqref{eq:MOL.05b} and
Appendix~\ref{sec:KT}.  But these terms, which originate from the
spatial discretization of a hyperbolic problem, have
a moderate CFL number, which scales linearly with the spatial stepsize,
so that the ``sub-problem''
$v_\K'=-H^0(v_\K)-H^1(v_\K)\mu$ is most efficiently implemented by an
explicit time marching
scheme.

To couple these contrary requirements, we introduce a $\halb$-explicit
IMEX-Runge-Kutta time-discretization for coupled DAEs of the form
\eqref{eq:HDAE}.  For $\halb$-explicit Runge-Kutta schemes for DAE
problems, see \cite{HairerLubichRoche:1989}.

To simplify the notation we write $V$ for $v_\K$ and restrict
the discussion mainly to \eqref{eq:HDAE1} and
\eqref{eq:HDAE2}, because \eqref{eq:HDAE3} and
\eqref{eq:HDAE4} decouple.  Hence, consider an ordinary DAE of the form
\begin{equation}\label{eq:DAEOrth}
  \begin{aligned}
    V'&=P(V)-H^1(V)\mu-H^0(V)+G(V),\\
    0&=H^1(V)^\top \Bigl(P(V)-H^1(V)\mu-H^0(V)+G(V)\Bigr)
  \end{aligned}
\end{equation}
in the case of the orthogonal phase condition \eqref{eq:ophase}
(resp.~\eqref{eq:MOL.AlgO}) and
\begin{equation}\label{eq:DAEFix}
  \begin{aligned}
    V'&=P(V)-H^1(V)\mu-H^0(V)+G(V),\\
    0&=\Psi^\top V-\Psi^\top \wh{U}
  \end{aligned}
\end{equation}
in the case of the fixed phase condition \eqref{eq:fphase}
(resp.~\eqref{eq:MOL.AlgF}).

From now on, consider general DAEs of the form \eqref{eq:DAEOrth} or
\eqref{eq:DAEFix}, i.e.\ $P$, $H^1$, $H^0$, and $G$ are given nonlinear
functions of $V$, sufficiently smooth in their arguments.

Let two Butcher-tableaux be given,\begin{center}
Tableau 1: \begin{tabular}{c|c}
  $c$&$A$\\\hline\\[-.9em]&$b^\top$
\end{tabular}
\hspace*{3cm}
Tableau 2:
\begin{tabular}{c|c}
  $c$&$\wh{A}$\\\hline\\[-.9em]&$\wh{b}^\top$
\end{tabular}
\end{center}
with $c=(c_0,\dots,c_s)^\top$,
$A=\bigl(a_{ij}\bigr)_{i,j=0,\dots,s}$,
$\wh{A}=\bigl(\wh{a}_{ij}\bigr)_{i,j=0,\dots,s}$ and
$b=(a_{s0},\dots,a_{ss})^\top$,
$\wh{b}=(\wh{a}_{s0},\dots,\wh{a}_{ss})^\top$.  We assume, that the
Runge-Kutta method with Tableau~1 is explicit, and
the Runge-Kutta method with Tableau~2 is
diagonally implicit, i.e.\ $a_{ij}=0$ for all $j\ge i$ and
$\wh{a}_{ij}=0$ for all $j> i$.  

Now assume that at some time instance $\tau_n$ a consistent
approximation $V^n$ and $\mu^n$ of $V$ and $\mu$ is given. 
A step for \eqref{eq:DAEOrth} from $\tau_n$ to $\tau_{n+1}=\tau_n+h$ with
stepsize $h$ is then performed as follows:  
\[\text{Set}\quad V_0=V^n\]
and for $i=1,\dots,s$ the internal values $V_i$ and $\mu_{i-1}$ of the
scheme are given as solutions to
\begin{equation}\label{eq:SysOrtho}
  \left\{
    \begin{aligned}
      V_i&=V_0-h \sum_{\nu=0}^{i-1}a_{i\nu}\Bigl(H^{1}(V_\nu)\mu_\nu+
      H^0(V_\nu)-G(V_\nu)\Bigr) + h\sum_{\nu=0}^i \wh{a}_{i\nu}
      P(V_\nu),\\
      0&=H^1(V_{i-1})^\top\Bigl(P(V_{i-1})
      -H^1(V_{i-1})\mu_{i-1}-H^0(V_{i-1})+
      G(V_{i-1})\Bigr),
    \end{aligned}
  \right.\quad i=1,\dots,s.
\end{equation}
As approximation of $V$ at the new time instance $\tau_{n+1}$
\[\text{set}\quad V^{n+1}:=V_s.\]
A suitable approximation of $\mu^{n+1}$ at the new time instance is actually the
internal value $\mu_0$ of the next time instance, i.e.
$\mu=\mu^{n+1}$ solves $0=H^1(V^{n+1})^\top\bigl(P(V^{n+1})
-H^1(V^{n+1})\mu-H^0(V^{n+1})+
G(V^{n+1})\bigr)$.

Similarly, a step for \eqref{eq:DAEFix} is performed as follows:
\[\text{Set}\quad V_{0}=V^n\]
and for $i=1,\dots,s$ let the internal values $V_{i}$ and $\mu_{i-1}$ be given
as solutions to
\begin{equation}\label{eq:SysFix}
  \left\{
    \begin{aligned}
      V_{i}&=V_{0}-h \sum_{\nu=0}^{i-1}a_{i\nu}\Bigl(H^{1}(V_{\nu})\mu_{\nu}+
      H^0(V_{\nu})-G(V_{\nu})\Bigr) + h\sum_{\nu=0}^i \wh{a}_{i\nu}
      P(V_{\nu}),\\
      0&=\Psi^{\top}V_{i}-\Psi^\top \wh{U},
    \end{aligned}
  \right.\quad i=1,\dots,s.
\end{equation}
In this case we set
\begin{center}
  $V^{n+1}:=V_{s}$ and $\mu^{n+1}:=\mu_{s-1}$
\end{center}
as approximations of $V$ and $\mu$ at the new time instance
$\tau_{n+1}=\tau_n+h$.

\begin{remark}\label{rem:schemes}
  \begin{enumerate}[label=\textup{(\roman*)}]
    \item
      A suitable time-discretization for the ``parabolic sub-problem''
      $v_\K'=P(v_\K)$ is the Crank-Nicolson method,
      for which there is no CFL-restriction and the resulting full
      discretization of $v'=\laplace v$ is of second order.
    \item
      Concerning the semi-discrete ``hyperbolic sub-problem''
      $v_\K'=-H^0(v_\K)-H^1(v_\K)\mu$, it is shown in
      \cite[Thm.~5.1]{KurganovTadmor:2000} that an explicit
      Euler-Discretization with a suitable CFL-condition satisfies a
      stability property in form of a maximum principle
      (cf.~non-oscillatory).  It is possible to retain this property for
      higher order methods, if they can be written as convex
      combinations of explicit Euler steps.  A simple second order
      method, for which this is possible is Heun's method, see
      \cite[Cor.~5.1]{KurganovTadmor:2000}.  It was observed already in
      \cite{ShuOsher:1988} that Heun's method is the optimal (concerning
      the CFL-restrictions) second order explicit Runge-Kutta type
      scheme.
  \end{enumerate}
\end{remark}

For concreteness and motivated by Remark~\ref{rem:schemes}, we choose
\begin{equation}\label{eq:ButcherTab}
  c=\vect{0\\1\\1},\;A=\vect{0&0&0\\1&0&0\\\halb&\halb&0},
  b^\top=\bigl(\halb,\halb,0\bigr),\;
  \wh{A}=\vect{0&0&0\\\halb&\halb&0\\\halb&0&\halb},\;
  \wh{b}^\top=\bigl(\halb,0,\halb\bigr),
\end{equation}
i.e.\ we couple Heun's method with the Crank-Nicolson method.

As noted above, the ODEs \eqref{eq:HDAE3} and \eqref{eq:HDAE4} can be
solved in a post-processing step, but if these values are also required,
it is more convenient to do the calculation in parallel and use the same
explicit scheme with Tableau~1 because the intermediate stages
of $V$ and $\mu$ are already calculated.  Therefore, given the state
$g^n$ and $t^n$ at $\tau_n$, we perform a time step for \eqref{eq:HDAE3}
and \eqref{eq:HDAE4} from $\tau_n$ to $\tau_{n+1}=\tau_n+h$ by
\[\text{let}\quad g_{0}=g^n,\;t_0=t^n,\]
compute for $i=1,\ldots,s$
\[ g_i=g_{0}+h
    \sum_{\nu=0}^{i-1}a_{i\nu}r_\text{alg}(g_\nu,\mu_\nu),\quad
    t_i=t_{0}+h \sum_{\nu=0}^{i-1}a_{i\nu}r_\text{time}(g_{\nu})
\]
and set
\[
  g^{n+1}:=g_s,\quad t^{n+1}:=t_s.
\]
  
\subsection{Order of the Time Discretization}
Now we consider the local truncation error of our IMEX-RK scheme for DAEs
\eqref{eq:SysOrtho}, resp.~\eqref{eq:SysFix}, with tableaux
\eqref{eq:ButcherTab}.  Hence, consider 
an ordinary differential equation
\begin{equation}\label{eq:DAE.ODE}
    V'=P(V)+H(V,\mu),\quad
    g'=r_\text{alg}(\mu,g),\quad
    t'=r_\text{time}(g),
\end{equation}
coupled to a system of algebraic equations of the form 
\begin{subequations}\label{eq:DAE.phase}
  \begin{equation}\label{eq:DAE.ind1}
    0=\Psi(V,\mu),
  \end{equation}
  or of the form
  \begin{equation}\label{eq:DAE.ind2}
    0=\Psi(V).
  \end{equation}
\end{subequations}
We assume that for consistent initial data $V(0)=V^0$, $g(0)=g^0$,
$t(0)=t^0$, $\mu(0)=\mu^0$
a smooth solution $V\in\CC^3([0,T);\R^m)$,
$\mu\in\CC^3([0,T);\R^p)$ of (\ref{eq:DAE.ODE},\ref{eq:DAE.ind1}) resp.
(\ref{eq:DAE.ODE},\ref{eq:DAE.ind2}) exists.  Furthermore we assume
\begin{assumption}\label{ass:DAE.ind}
  \begin{enumerate}[label=\textup{(\roman*)}]
    \item\label{itm:DAE.ind1}
      In the case of \eqref{eq:DAE.ODE} with \eqref{eq:DAE.ind1} the matrix 
      $\partial_\mu \Psi\bigl(V(\tau),\mu(\tau)\bigr)$ is invertible for
      any $\tau\in[0,T)$.
    \item\label{itm:DAE.ind2}
      In the case of \eqref{eq:DAE.ODE} with \eqref{eq:DAE.ind2} the matrix 
      $\partial_V \Psi\bigl(V(\tau)\bigr)\partial_\mu
      H\bigl(V(\tau),\mu(\tau)\bigr)$ is invertible for any
      $\tau\in[0,T)$.
  \end{enumerate}
\end{assumption}
It is easy to check that under Hypothesis~\ref{ass:DAE.ind} the system
(\ref{eq:DAE.ODE}, \ref{eq:DAE.ind1}) is a DAE of (differentiation) index $1$
and the system
(\ref{eq:DAE.ODE}, \ref{eq:DAE.ind2}) is a DAE of (differentiation) index $2$
(cf. \cite[Ch.~VII]{HairerWanner:1996}).

Given consistent data $V^n, g^n, t^n$ of the DAE (\ref{eq:DAE.ODE},
\ref{eq:DAE.ind1}) or (\ref{eq:DAE.ODE}, \ref{eq:DAE.ind2})
at some time-instance $\tau_n$, a step of size $h$ of the method with
these data takes the following explicit form:
Solve the system
\begin{subequations}\label{eq:step}
  \begin{equation}\label{eq:step.init}
    V_0=V^n,\quad g_0= g^n,\quad t_0=t^n,
  \end{equation}
  \begin{align}
    &\left[
      \begin{aligned}
        V_1&=V_0+\tfrac{h}{2} P(V_0)+\tfrac{h}{2} P(V_1)+h H(V_0,\mu_0),\\
        0&=\begin{cases} \Psi(V_0,\mu_0),\quad&\text{case
            \eqref{eq:DAE.ind1}, or}\\
            \Psi(V_1),\quad&\text{case \eqref{eq:DAE.ind2}},
        \end{cases}\\
        g_1&=g_0+h r_\text{alg}(\mu_0,g_0),\\
        t_1&=t_0+h r_\text{time}(g_0),
      \end{aligned}
    \right.\\
    &\left[
      \begin{aligned}
        V_2&=V_0+\tfrac{h}{2} P(V_0)+\tfrac{h}{2} P(V_2)+
        \tfrac{h}{2} H(V_0,\mu_0)+\tfrac{h}{2} H(V_1,\mu_1),\\
        0&=\begin{cases} \Psi(V_1,\mu_1),\quad&\text{case
            \eqref{eq:DAE.ind1}, or}\\
            \Psi(V_2),\quad&\text{case \eqref{eq:DAE.ind2}},
        \end{cases}\\
        g_2&=g_0+\tfrac{h}{2} r_\text{alg}(\mu_0,g_0)+
        \tfrac{h}{2} r_\text{alg}(\mu_1,g_1),\\
        t_2&=t_0+\tfrac{h}{2} r_\text{time}(g_0)+
        \tfrac{h}{2} r_\text{time}(g_1).
      \end{aligned}
    \right.
  \end{align}
\end{subequations}
Finally, the values at the new time-instance $\tau^{n+1}=\tau^n+h$ are
given by
\begin{equation}\label{eq:stepfin}
  V^{n+1}:=V_2,\; g^{n+1}:=g_2,\;t^{n+1}:=t_2,\;\text{and}
  \mu^{n+1}\; \begin{cases}\text{ by solving }
    \Psi(V^{n+1},\mu^{n+1})=0,&\text{case \eqref{eq:DAE.ind1}},\\
    \mu_1,&\text{case \eqref{eq:DAE.ind2}.}\end{cases}
\end{equation}
Since the $g$- and $t$-equation decouple from the system, we first
consider the $V$ and $\mu$ variables separately.  To analyze the local
error, we assume that $(V_\star,\mu_\star)$ is a given consistent
value at $\tau=0$, so that the DAE satisfies all assumptions
from above with $(V^0,\mu^0)$
replaced by $(V_\star,\mu_\star)$.  For brevity,
a subindex $\star$ denotes the evaluation of a function at $\tau=0$,
e.g.\ $(V_\star,\mu_\star)=(V(0),\mu(0))$, $\partial_V
P_\star
  =\tfrac{\partial}{\partial V}P(V(0))$.

\textbf{Taylor expansion of the exact solution.}
First we consider the differential variables $V$ and obtain from the ODE
\eqref{eq:DAE.ODE} and anticipating $\mu(0)=\mu_\star$
(see \eqref{eq:mu_exp1a}, resp.\ \eqref{eq:mu_exp2b}, below)
\begin{subequations}\label{eq:V_exp}
  \begin{align}
    V(0)&=V_\star,\label{eq:V_expa}\\
    V'(0)&=P_\star+H_\star=:V'_\star,\label{eq:V_expb}\\
    V''(0)&=\partial_V P_\star V'_\star+\partial_V H_\star V'_\star+
    \partial_\mu H_\star \mu'(0)=: V''_\star.\label{eq:V_expc}
  \end{align}
\end{subequations}
Similarly, we can use the differential equation \eqref{eq:DAE.ODE} together
with the algebraic constraint \eqref{eq:DAE.ind1} to obtain in the index-$1$
case for the algebraic variable
\begin{subequations}\label{eq:mu_exp1}
  \begin{align}
    0&=\Psi(V_\star,\mu_{\star}),\quad\text{(locally unique solvable for
      $\mu_\star$ by Hypothesis~\ref{ass:DAE.ind}~\ref{itm:DAE.ind1})},\label{eq:mu_exp1a}\\
    \mu'(0)&=-\bigl(\partial_\mu \Psi_\star\bigr)^{-1} \partial_V \Psi_\star
    \bigl(P_\star+H_\star\bigr)=: \mu'_{\star},\label{eq:mu_exp1b}\\
    \mu''(0)&=-\bigl(\partial_\mu \Psi_\star\bigr)^{-1}
    \Bigl\{\partial_V^2\Psi_\star {V'_\star}^2+\partial_V \Psi_\star V''_\star
      + 2\partial_V\partial_\mu \Psi_\star V'_\star \mu'_\star+\partial_\mu^2
      \Psi_\star {\mu'_\star}^2\Bigr\}=:\mu''_{\star}.
  \end{align}
\end{subequations}
For the index-$2$ case (i.e.~\eqref{eq:DAE.ind2}) the first coefficients of
the Taylor expansion of the algebraic variables $\mu$ are obtained by
differentiating the algebraic condition \eqref{eq:DAE.ind2} and using the
ODE \eqref{eq:DAE.ODE} to find
\begin{subequations}\label{eq:mu_exp2}
  \begin{align}
    0&=\Psi(V_\star),\\
    0&=\partial_V \Psi_\star (P_\star+H_\star),\quad\text{(locally
      unique solvable for $\mu_\star$ by
      Hypothesis~\ref{ass:DAE.ind}~\ref{itm:DAE.ind2})},\label{eq:mu_exp2b}\\
    \mu'(0)&= -\bigl(\partial_V \Psi_\star \partial_\mu H_\star\bigr)^{-1}
    \Bigl\{\partial_V^2 \Psi_\star {V'_\star}^2+\partial_V \Psi_\star
      \partial_V P_\star V'_\star+\partial_V \Psi_\star \partial_V H_\star
    V'_\star\Bigr\}=:\mu'_\star.\label{eq:mu_exp2c}
  \end{align}
\end{subequations}

\textbf{Taylor expansion of the numerical solution.}  We assume that the
numerical solution and all intermediate stages depend smoothly on the
step-size $h$.  To emphasize this dependence, we explicitly include
the dependence on $h$ in the notation.

First consider the differential variable $V$,
anticipating that the algebraic variables $\mu_0(h)$ and $\mu_1(h)$ are
known and satisfy $\mu_0(0)=\mu_1(0)=\mu_\star$, which will be justified
in \eqref{eq:mu001} and \eqref{eq:mu101}, resp.\ \eqref{eq:mu002}
and \eqref{eq:mu102}.
By \eqref{eq:step} the numerical values $V_1(h)$ and $V_2(h)$ satisfy
\begin{equation}\label{eq:Vdisc}
  \begin{aligned}
    V_1(h)&=V_0+\tfrac{h}{2}P(V_0)+\tfrac{h}{2}P(V_1(h))+h
    H(V_0,\mu_0(h))\\
  &=V_\star+\tfrac{h}{2} P_\star+\tfrac{h}{2}P(V_1(h))+h
  H(V_\star,\mu_0(h)),\\
  V_2(h)&=V_0+\tfrac{h}{2}P(V_0)+\tfrac{h}{2}P(V_2(h))+
  \tfrac{h}{2}H(V_0,\mu_0(h))+\tfrac{h}{2}H(V_1(h),\mu_1(h))\\
  &=V_\star+\tfrac{h}{2}P_\star+\tfrac{h}{2}P(V_2(h))
  +\tfrac{h}{2}H(V_\star,\mu_0(h))+\tfrac{h}{2} H(V_1(h),\mu_1(h)).
\end{aligned}
\end{equation}
For the differential variables we then obtain for the first intermediate
stage:
\begin{subequations}\label{eq:numV1_exp}
  \begin{align}
    V_1(0)&=V_\star,\\
    V'_1(h)&=\tfrac{1}{2} P_\star+\tfrac{1}{2} P(V_1(h))
    +\tfrac{h}{2} \partial_V P(V_1(h))
    V_1'(h)+H(V_\star,\mu_0(h))+h\partial_\mu H(V_\star,\mu_0(h))
    \mu_0'(h),\\
    V'_1(0)&=P_\star+H_\star=V'_\star,\\
    V''_1(0)&=\partial_V P_\star (P_\star+H_\star)+2\partial_\mu H_\star
    \mu_0'(0).
  \end{align}
\end{subequations}
Similarly, for the final (second) stage we obtain
\begin{subequations}\label{eq:numV2_exp}
  \begin{align}
    V_2(0)&=V_\star,\\
    V'_2(h)&=\tfrac{1}{2} P_\star+\tfrac{1}{2} P(V_2(h))
    +\tfrac{h}{2} \partial_V P(V_2(h))
    V_2'(h)+\tfrac{1}{2}H(V_\star,\mu_0(h))+\tfrac{h}{2} 
    \partial_\mu H(V_\star,\mu_0(h))\mu_0'(h)\\
    &\qquad+ \tfrac{1}{2}H(V_1(h),\mu_1(h))
    +\tfrac{h}{2} \partial_V H(V_1(h),\mu_1(h))
    V_1'(h)+\tfrac{h}{2}\partial_\mu H(V_1(h),\mu_1(h)) \mu_1'(h),\notag\\
    V'_2(0)&=P_\star+H_\star=V'_\star,\\
    V''_2(0)&=\partial_V P_\star (P_\star+H_\star)+\partial_V
    H_\star(P_\star+H_\star)+\partial_\mu
    H_\star(\mu_0'(0)+\mu_1'(0)).\label{eq:numV2_expd}
  \end{align}
\end{subequations}

Now we consider the Taylor expansion of the numerical solution of the
algebraic variables.  We begin with the \emph{index-$1$ case},
i.e.~\eqref{eq:Vdisc} is closed by imposing the algebraic constraints
\begin{equation}\label{eq:algdisc1}
  0=\Psi(V_0(h),\mu_0(h)),\quad 0=\Psi(V_1(h),\mu_1(h)).
\end{equation}
Since $V_0(h)=V_\star$ for all $h\ge 0$, \eqref{eq:algdisc1} implies 
because of Hypothesis~\ref{ass:DAE.ind}~\ref{itm:DAE.ind1}
\begin{equation}\label{eq:mu001}
  \mu_0(h)=\mu_\star\text{ and }
  \mu_0'(h)=0\;\forall h\ge 0.
\end{equation}
This justifies \eqref{eq:numV1_exp} and shows
$V_1''(0)=\partial_V P_\star (P_\star+H_\star)$.  Then we can
calculate the Taylor expansion of
$\mu_1(h)$ around $0$ from \eqref{eq:algdisc1} to find
\begin{equation}\label{eq:mu101}
  \begin{aligned}
    0&=\Psi(V_1(h),\mu_1(h))&\;\text{which implies}\;
    \mu_1(0)&=\mu_\star\;\text{since $V_1(0)=V_\star$},\\
    0&=\partial_V\Psi_1 V_1'+\partial_\mu \Psi_1 \mu_1'&\;\text{which implies}\;
    \mu_1'(0)&=-\bigl(\partial_\mu \Psi_\star\bigr)^{-1} \partial_V
    \Psi_\star (P_\star+H_\star)=\mu_\star'\;\text{from
      \eqref{eq:mu_exp1b}}.
  \end{aligned}
\end{equation}
Thus, inserting the findings for $\mu_0'(0)$ and $\mu_1'(0)$
into the Taylor expansion of $V_2$,
\eqref{eq:numV2_exp}, shows
$V_2(0)=V_\star$, $V_2'(0)=V_\star'$, and $V_2''(0)=V_\star''$,
which proves for the differential variable $V_2$ second order
consistency,
\begin{equation}\label{eq:Verr}
  V_2(h)=V(h)+{\mathcal O}(h^3)\quad\text{as $h\searrow 0$}.
\end{equation}
Setting $\mu_2(h)$ as solution of $\Psi(V_2(h),\mu_2(h))=0$ implies
also second order for the algebraic variable
\begin{equation}\label{eq:muerr_ind1}
  \mu_2(h)=\mu(h)+{\mathcal O}(h^3)\quad\text{as $h\searrow 0$}.
\end{equation}

In the \emph{index-$2$ case} \eqref{eq:Vdisc} is closed by appending the
algebraic constraints
\begin{equation}\label{eq:algdisc2}
  0=\Psi\bigl(V_1(h)\bigr),\quad 0=\Psi\bigl(V_2(h)\bigr).
\end{equation}
Differentiating the first of these two equations yields
$0=\partial_V \Psi(V_1(h)) V_1'(h)$ and at $h=0$
\begin{equation}\label{eq:mu002}
  0=\partial_V \Psi_\star \bigl(P_\star+H(V_\star,\mu_0(0))\bigr),
\end{equation}
which by Hypothesis~\ref{ass:DAE.ind}~\ref{itm:DAE.ind2} has the locally
unique solution
$\mu_0(0)=\mu_\star$.  Considering the second derivative at
$h=0$ leads to
\[0=\partial_V^2 \Psi_\star(P_\star+H_\star)^2+\partial_V \Psi_\star
  \bigl(\partial_V P_\star (P_\star+H_\star)+2\partial_\mu H_\star
\mu_0'(0)\bigr),\]
so that 
\begin{equation}
  \label{eq:mu0'0equiv}
  \mu_0'(0)=-\bigl(\partial_V\Psi_\star\partial_\mu H_\star\bigr)^{-1}
  \Bigl\{ \tfrac{1}{2}\partial_V^2
  \Psi_\star(P_\star+H_\star)^2+\tfrac{1}{2}\partial_V
  \Psi_\star \bigl(\partial_V P_\star(P_\star+H_\star)\bigr)\Bigr\}.
\end{equation}
Similarly, differentiating the second equation in \eqref{eq:algdisc2}
once and evaluating at $h=0$ yields
\[0=\partial_V\Psi_\star \bigl(P_\star+\tfrac{1}{2}H_\star
+\tfrac{1}{2}H(V_\star,\mu_1(0))\bigr),\]
which has the locally unique solution
\begin{equation}\label{eq:mu102}
  \mu_1(0)=\mu_\star.
\end{equation}
Considering the second derivative at $h=0$, inserting
\eqref{eq:numV2_expd} and comparing with \eqref{eq:mu_exp2c} shows
\begin{equation}
    \mu_0'(0)+\mu_1'(0)=\mu'_\star.
\end{equation}
Therefore, we obtain
again \eqref{eq:Verr} for $V_2$.
But for the algebraic variables we only find the estimate
\begin{equation}\label{eq:muerr_ind2}
  \mu_1(h)+{\mathcal O}(h)=\mu_2(h)+{\mathcal
  O}(h)=\mu(h)\quad\text{as\;}h\searrow 0,
\end{equation}
which is a severe order reduction.

Nevertheless, when the group variables $g$ and the transformed time
$t$ are calculated in parallel with $V$ and $\mu$, these variables are
again second order accurate, i.e.
\begin{equation}\label{eq:othererr}
  g_2(h)=g(h)+{\mathcal O}(h^3),\quad t_2(h)=t(h)+{\mathcal
  O}(h^3)\quad\text{as $h\searrow 0$}.
\end{equation}
This can be seen, by either including the equations
$g'=r_\text{alg}(\mu,g)$ and $t'=r_\text{time}(g)$ into the
$V$-equation, or by performing a similar analysis as above.

Summarizing the above analysis we obtain.
\begin{theorem}\label{thm:2ndorder}
  Assume that the DAE consisting of the ODE system \eqref{eq:DAE.ODE} and
  either \eqref{eq:DAE.ind1} or \eqref{eq:DAE.ind2} with consistent initial
  data $(V_\star,\mu_\star,g_\star,t_\star)$ at $\tau=\tau_n$ 
  has a smooth solution in
  some non-empty interval $[\tau_n,T)$ and satisfies
  Hypothesis~\ref{ass:DAE.ind}.  Then the method \eqref{eq:step},
  \eqref{eq:stepfin} is
  consistent of second order at the differential variables, i.e.
  \eqref{eq:Verr} and \eqref{eq:othererr} hold.
  Moreover, in the case of \eqref{eq:DAE.ind1} also the algebraic
  variables $\mu$ are
  second order consistently approximated, i.e.
  \eqref{eq:muerr_ind1} holds.
\end{theorem}
\begin{remark}
  We observe here the well-known difficulty with higher index problems, that
  one faces a loss in the order of the approximation of the algebraic
  variables.  But note, that one is often mainly interested in the
  behavior of the original solution, so that the approximation of $\mu$ is
  not important but the approximations of $V$, $g$, and $t$ are crucial.
  
  Moreover, when the main focus is on the asymptotic behavior
  and the final rest state, one is interested in the limit
  $\lim_{\tau\to\infty}\mu(\tau)$ but
  not on its actual evolution.  The approximation of this is actually not
  dependent on the time-discretization but only on the spatial
  discretization, because we are using a Runge-Kutta-type method.
  In fact, in this case it might even be better to use a scheme which
  has a less restrictive CFL condition but maybe has a lower order in
  the time approximation.  See \cite{Shu:1988} for such ideas in the
  case of hyperbolic conservation laws.
\end{remark}
  
\subsection{Efficiently solving the Runge-Kutta Equations}
We remark that in Burgers' case and 
similar cases with a linear operator $P$, the actual equations
\eqref{eq:SysOrtho} and \eqref{eq:SysFix} from the
IMEX-time-discretization can be solved very efficiently.  To see this,
we write
\eqref{eq:SysOrtho} in block-matrix form
\begin{equation}\label{eq:FactorOrtho}
  \begin{pmatrix}I-h\,\wh{a}_{i,i} P&h\,a_{i,i-1} \,H^1(V_{i-1})\\
    0&H^1(V_{i-1})^\top H^1(V_{i-1})
  \end{pmatrix} \begin{pmatrix} V_{i}\\\mu_{i-1}\end{pmatrix}
  = \begin{pmatrix} R^1_i\\
    R^2_i
  \end{pmatrix},\quad i=1,\dots,s,
\end{equation}
where $R^1_i=R^1_i(V_{0},\dots,V_{i-1},\mu_{0},\dots,\mu_{i-2})$ and
$R^2_i=R^2_i(V_{i-1})$ are given by
\begin{align*}
  R^1_i&=V_{0}-h\,\Bigl(\sum_{\nu=0}^{i-2} a_{i\nu}
  H^1(V_{\nu})\mu_{\nu}+\sum_{\nu=0}^{i-1} a_{i\nu}
  \bigl(H^0(V_{\nu})-G(V_{\nu})\bigr)- \sum_{\nu=0}^{i-1}
  \wh{a}_{i\nu} P(V_{\nu})\Bigr),\\
  R^2_i&=H^1(V_{i-1})^\top\bigl(P(V_{i-1})-
  H^0(V_{i-1})+G(V_{i-1})\bigr).
\end{align*}
A solution to \eqref{eq:FactorOrtho} can be calculated by solving for
each $i=1,\dots,s$ one small
linear system with the $(d+1)\times(d+1)$-matrix $H^1(V^{(i-1)})^\top
H^1(V^{(i-1)})$ and one large linear
system with the matrix $I-h\, \wh{a}_{i,i}P$.

Similarly, \eqref{eq:SysFix} can be written in the block-matrix form
\begin{equation}\label{eq:FactorFix}
  \begin{pmatrix}I-h\, \wh{a}_{i,i}\,P&h\,a_{i,i-1}
    \,H^1(V_{i-1})\\
    \Psi^\top&0
  \end{pmatrix} \begin{pmatrix} V_{i}\\\mu_{i-1}\end{pmatrix}
  = \begin{pmatrix} R^1_i\\
    R^2
  \end{pmatrix},\quad i=1,\dots,s,
\end{equation}
where $R^1_i$ is given as above and $R^2=\Psi^\top \wh{U}$ does not
depend on the stage $i$.  It is possible to obtain a solution to
\eqref{eq:FactorFix} by solving $d+2$ large
linear systems with $I-h\,\wh{a}_{i,i}P$ and one small linear
$(d+1)\times(d+1)$ system:
\begin{equation*}
  \begin{aligned}
    A_\star &= \bigl(I-h\,\wh{a}_{i,i}P\bigr)^{-1} h
    a_{i,i-1}H^1(V_{i-1}),&&\text{($(d+1)\times$ large)}\\
    V_\star &=  \bigl(I-h\,\wh{a}_{i,i}P\bigr)^{-1}
    R_i^1,&&\text{($1\times$ large)}\\
    \mu_{i-1}&=\bigl(\Psi^\top A_\star\bigr)^{-1} \bigl(
    \Psi^\top V_\star-R^2\bigr),&&\text{($1\times$ small)}\\
    V_i&=V_\star-A_\star\mu_{i-1}.
  \end{aligned}
\end{equation*}
Therefore, if $h$ does not vary
during computation, factorizations of $I-h\,\wh{a}_{i,i}P$,
$i=1,\dots,s$ can be
calculated in a preprocessing step and then the subsequent time steps
are rather cheap.
Note that for our specific scheme with Butcher
tableaux \eqref{eq:ButcherTab} holds $s=2$ and
$\wh{a}_{1,1}=\wh{a}_{2,2}=\tfrac{1}{2}$, so that only one matrix
factorization is needed, namely of $I-h \tfrac{1}{2} P$.


\section{Numerical experiments}\label{sec:5}
In the current paper, we are primarily interested in the actual
discretization of the freezing PDAE \eqref{eq:CoCauchy} and the error
introduced by the discretization.  We present results of several numerical
experiments which support second order convergence of the method for the
full discretization of the PDAE for our scheme.

For further numerical findings we refer to \cite{Rottmann:2016a},
where we use the same scheme and see that our method enables us to do
long time simulations for Burgers' equation on fixed (bounded)
computational domains also for different parameter values of $p>1$ in
\eqref{eq:burgersdd}.  We also note that our method is also able to
capture the meta-stable solution behavior present in
\eqref{eq:burgersdd} for the case $0<\nu<<|a|$ and again refer to
\cite{Rottmann:2016a}.  Note that in that paper we completely ignored
the errors introduced by the discretization.

In the following, we restrict to the special choice $p=\tfrac{d+1}{d}$
in \eqref{eq:burgersdd}.  For convenience we explicitly state the
systems for freezing similarity solutions of the ``\emph{conservative}''
Burgers' equation (i.e.\ $p=\tfrac{d+1}{d}$) in one and two
spatial dimensions.  These are the systems we numerically solve with the
scheme proposed in Sections~\ref{sec:3}--\ref{sec:4}.
In the \textbf{1d case} we solve the PDAE system
\begin{equation}\label{eq:1dexpl}
  \left\{
  \begin{aligned}
    v_\tau&=\nu v_{\xi\xi}-\tfrac{1}{2}\left(v^2\right)_\xi +\mu_1 (\xi
    v)_\xi + \mu_2 v_\xi,&v(0)&=u_0,\\
    0&=\Psi(v,\mu)\in\R^2,\\
    \alpha_\tau&=\alpha \mu_1,\;
    b_\tau=\alpha \mu_2,\;
    t_\tau=\alpha^2,
    &\alpha(0)&=1,\,b(0)=0,\,t(0)=0.
  \end{aligned}
  \right.
\end{equation}
And the functional $\Psi$ is either  given
by 
\[
  \Psi(v,\mu)=\vect{\int_\R (\xi v)_\xi\cdot\Bigl( \nu
    v_{\xi\xi}-\tfrac{1}{2}(v^2)_\xi +\mu_1(\xi v)_\xi +\mu_2
    v_\xi\Bigr)\,d\xi\\
    \int_\R v_\xi\cdot\Bigl( \nu
    v_{\xi\xi}-\tfrac{1}{2}(v^2)_\xi +\mu_1(\xi v)_\xi +\mu_2
  v_\xi\Bigr)\,d\xi}\quad\text{(orthogonal phase condition)}
\]
or  by 
\[\Psi(v) = \vect{ \int_\R(\xi v)_\xi\cdot\bigl(\wh{u}-v\bigr)\,d\xi\\
\int_\R v_\xi\cdot\bigl(\wh{u}-v\bigr)\,d\xi}\quad \text{(fixed phase
condition)}.\]
Similarly, in the \textbf{2d case} we numerically solve
\begin{equation}\label{eq:2dexpl}
  \left\{\begin{aligned}
    v_\tau&=\underbrace{\nu \,\laplace
    v-\tfrac{2}{3}\left(\partial_{\xi_1}+\partial_{\xi_2}\right)|v|^\frac{3}{2}
    +\tfrac{1}{2}\mu_1 \bigl((\xi_1 v)_{\xi_1}+(\xi_2
    v)_{\xi_2}\bigr)+\mu_2 v_{\xi_1}+\mu_3
    v_{\xi_2}}_{=:\text{RHS}(v,\mu)},\quad v(0)=u_0\\
    0&=\Psi(v,\mu)\in\R^3,\\
    \alpha_\tau&=\alpha \mu_1,\;
    b_{1,\tau}=\alpha^\frac{1}{2}
    \mu_2,\;
    b_{2,\tau}=\alpha^\frac{1}{2}
    \mu_3,\;
    t_\tau=\alpha,\quad
    \alpha(0)=1,\,b_1(0)=0,\,b_2(0)=0,\,t(0)=0.
  \end{aligned}\right.
\end{equation}
The functional $\Psi$ is either given
by
\[
  \Psi(v,\mu)=\vect{\int_{\R^2} \bigl((\xi_1 v)_{\xi_1}+(\xi_2
  v)_{\xi_2}\bigr)\cdot\text{RHS}(v,\mu)\,d\xi\\
  \int_{\R^2} v_{\xi_1}\cdot\text{RHS}(v,\mu)\,d\xi\\
  \int_{\R^2} v_{\xi_2}\cdot\text{RHS}(v,\mu)\,d\xi}
  \quad\text{(orthogonal phase condition)}
\]
or by
\[\Psi(v) = \vect{ \int_{\R^2}\bigl((\xi_1 v)_{\xi_1}+(\xi_2
  v)_{\xi_2}\bigr)\cdot\bigl(\wh{u}-v\bigr)\,d\xi\\
\int_{\R^2} v_{\xi_1}\cdot\bigl(\wh{u}-v\bigr)\,d\xi\\
\int_{\R^2} v_{\xi_2}\cdot\bigl(\wh{u}-v\bigr)\,d\xi}\quad \text{(fixed
  phase condition)}.\]

In the following we write $\Delta \tau$ for the time step-size $h$ and
$\Delta \xi$, resp. $\Delta \xi_1$ and $\Delta \xi_2$, for the spatial
step-sizes.

\subsection{1d-Experiments}
We choose a fixed spatial domain $\Omega$, as explained in
Section~\ref{sec:3.2}.
Due to the nonlinearity, and because the algebraic variables
$\mu$ depend implicitly on the solution, we choose 
a very rough upper estimate for the value of $\mathbf{a}$ in
\eqref{eq:MOL.05} for simplicity.  We take (cf.~\eqref{eq:MOL.05e})
\begin{equation}\label{eq:wspeed1d}
  \mathbf{a}\le\sup_{\xi\in\Omega}|v(\xi)|+
  |\mu_1|\max\{|\xi|:\xi\in\Omega\}+|\mu_2|,
\end{equation}
where $\mu_1$, $\mu_2$ are the current approximations obtained in the
last time step, so that the number $\mathbf{a}$ is
updated after each time step.

As time step-size we choose $\Delta\tau$ which satisfies
\begin{equation}\label{eq:deltatstep}
  \Delta \tau\le \frac{\Delta \xi}{\mathbf{a}}\cdot\lambda_\mathrm{CFL},
\end{equation}
so the time step-size may increase or
decrease during the calculation, depending on the evolution of
$\mathbf{a}$ from \eqref{eq:wspeed1d}.  Experimentally, we found that 
$\lambda_\mathrm{CFL}=\tfrac{1}{3}$ is a suitable choice for all spatial
step-sizes $\Delta \xi$ and all viscosities $\nu\ge 0$.
In our code we actually do not update $\Delta \tau$ after each step, but
only if \eqref{eq:deltatstep} is violated or $\Delta \tau$ might be
enlarged significantly according to \eqref{eq:deltatstep}.  

\begin{figure}
  \centering
  \subfigure[$\nu=0.4$]{\includegraphics[height=4cm]{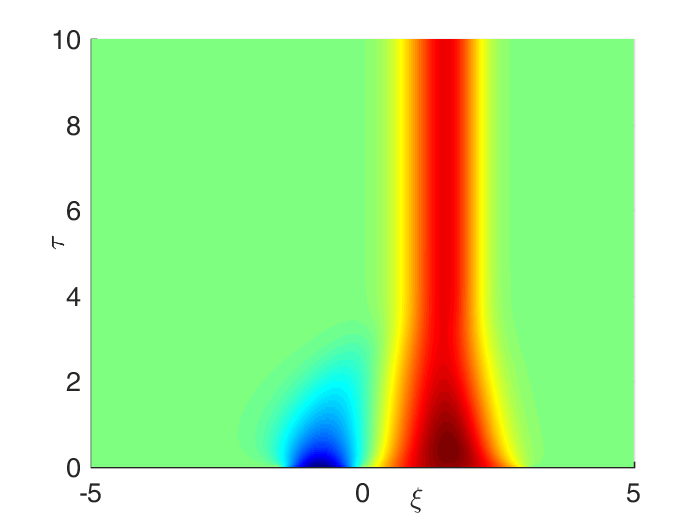}}%
  \subfigure[$\nu=0.01$]{\includegraphics[height=4cm]{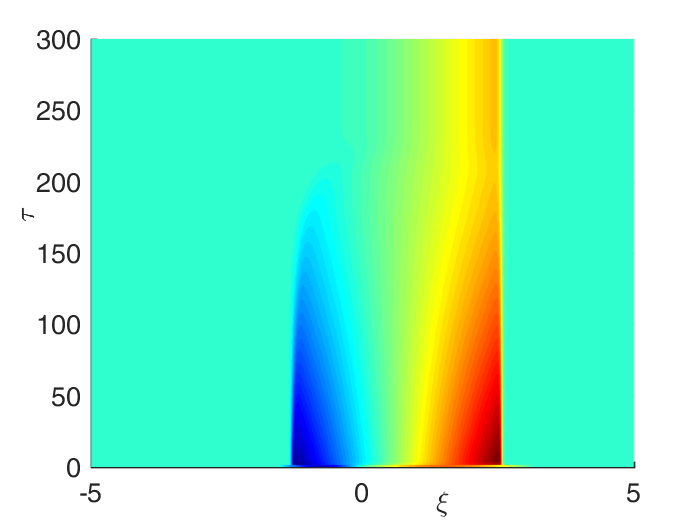}}%
  \subfigure[$\nu=0$]{\includegraphics[height=4cm]{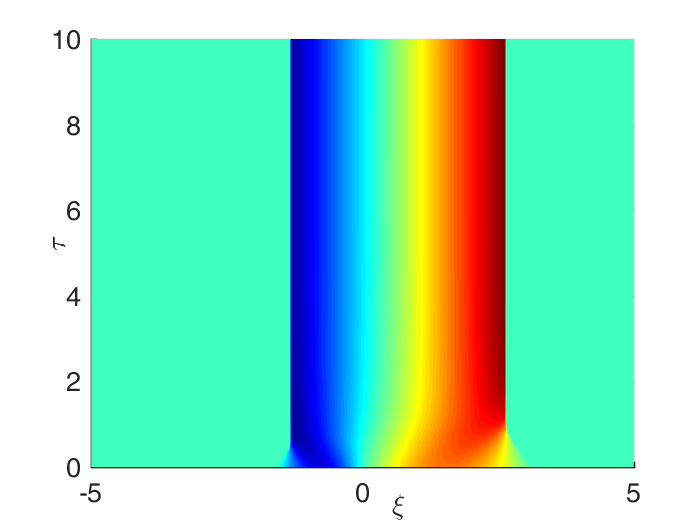}}
  \subfigure[$\nu=0.4,\; \tau=10, t\approx4.6\cdot10^{11}$]{\includegraphics[height=4cm]{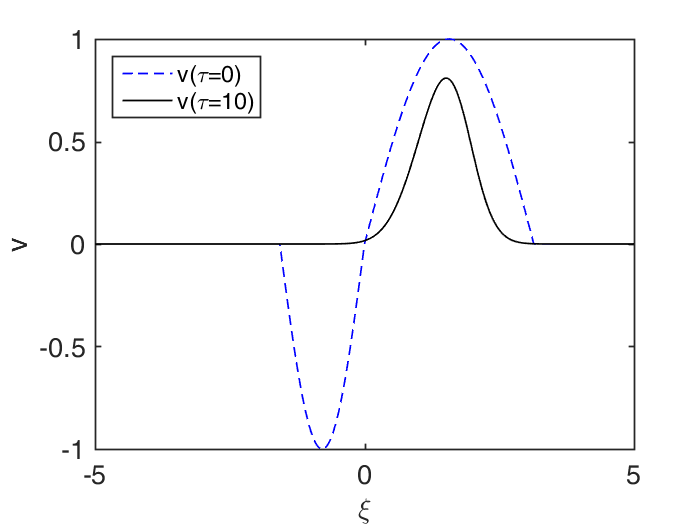}\label{fig:1dltp04}}%
  \subfigure[$\nu=0.01,\; \tau=300, t\approx5.22\cdot 10^{54}$]{\includegraphics[height=4cm]{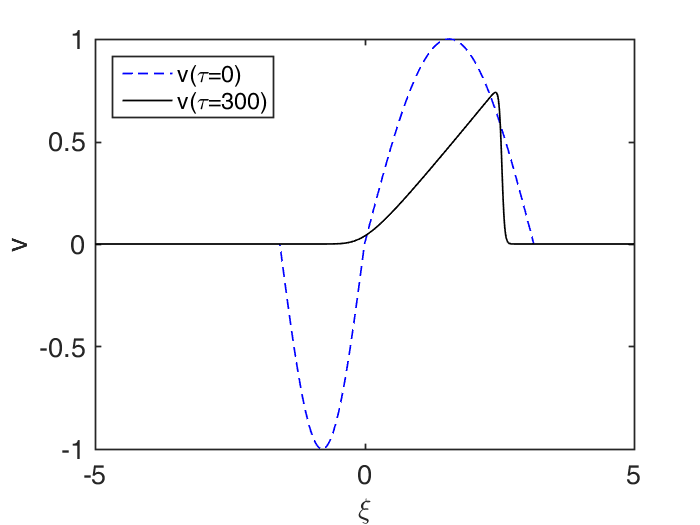}\label{fig:1dltp001}}%
  \subfigure[$\nu=0,\; \tau=10, t\approx1.75\cdot 10^3$]{\includegraphics[height=4cm] {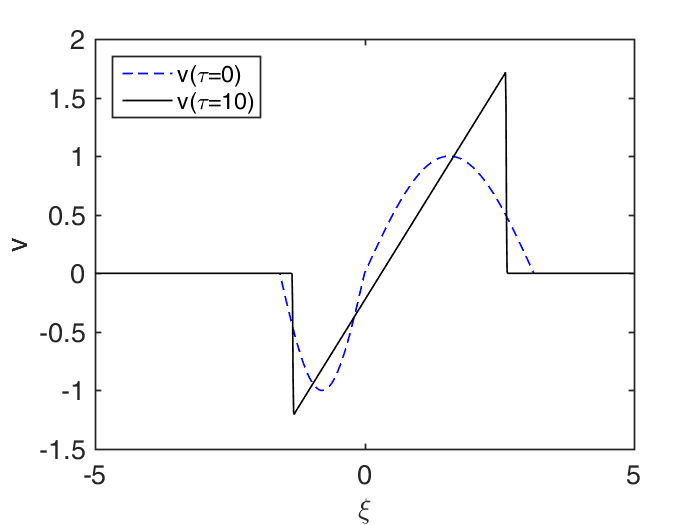}\label{fig:1dltp0}}
  \caption{Time evolution (a)--(c) and final states after the solution
    stabilized (d)--(f) for the freezing method for 1d-Burgers' equation
    with different sizes of viscosity.  The initial value is
    given by the dashed line, the final state by the solid line.  All
    plots are in the scaled and translated coordinates of the freezing
  method.}\label{fig:1dltplots}
\end{figure}
In Fig.~\ref{fig:1dltplots}(a)--(c) we show the time evolution of the
solution to the freezing PDAE \eqref{eq:1dexpl} for $\nu=0.4$,
$\nu=0.01$ and $\nu=0$.  The initial function is
$u_0(x)=\sin(2x)1_{[-\frac{\pi}{2},0]}+\sin(x)1_{[0,\pi]}$, where
$1_{A}$ is the indicator function.  One observes that the solution
converges to a stationary profile as time increases.
We stopped the calculation when the solution did not change anymore.
The initial value together with the final state
of these simulations is shown in Fig.~\ref{fig:1dltplots}(d)--(f).
Note that in all plots the solution is given in the co-moving
coordinates of \eqref{eq:1dexpl} and the solution in the original
coordinates is obtained by the transformation
$u(x,t)=\tfrac{1}{\alpha(\tau)}v(\tfrac{x-b(\tau)}{\alpha(\tau)},\tau)$.
At the final time the algebraic variables $\alpha$ and $\mu$ have the
values $\alpha(10)\approx1.2\cdot 10^6$ and $b(10)\approx-1.6\cdot 10^6$
in Fig.~\ref{fig:1dltp04}, $\alpha(300)\approx 1.3\cdot 10^{27}$ and
$b(300)\approx-1.8\cdot 10^{26}$ in Fig.~\ref{fig:1dltp001},
$\alpha(10)\approx 36$ and $b(10)\approx -11$ for Fig.~\ref{fig:1dltp0}.

\textbf{Order of the method.}
Now we numerically check the order of our method.  For this we calculate
the solution for different step-sizes $\Delta\xi$ until time $\tau=1$ and
compare the final state with a reference
solution, obtained by using a much smaller step-size $\Delta
\xi=0.0005$.  We calculate on the fixed spatial domain
$[R^-,R^+]$ with no-flux boundary conditions as described in
Section~\ref{sec:3}.  
\begin{figure}[htb]
  \centering
  \subfigure[$\nu=1$, orthogonal phase]{\includegraphics[height=5cm]%
    {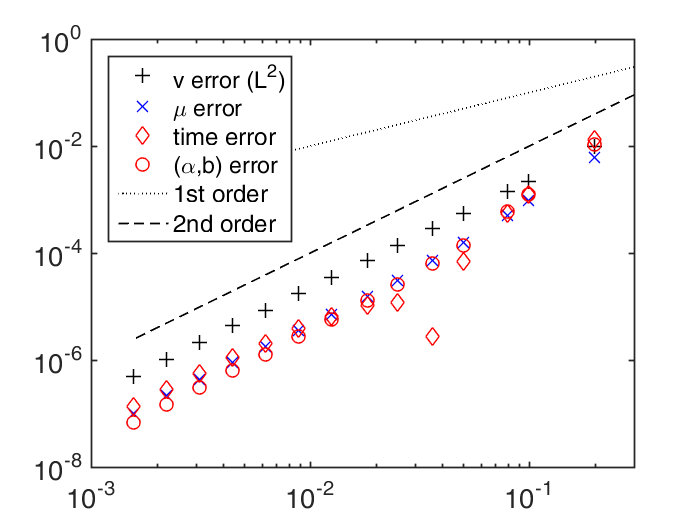}\label{fig:CLO}}
  \subfigure[$\nu=0.01$, orthogonal phase]{\includegraphics[height=5cm]%
    {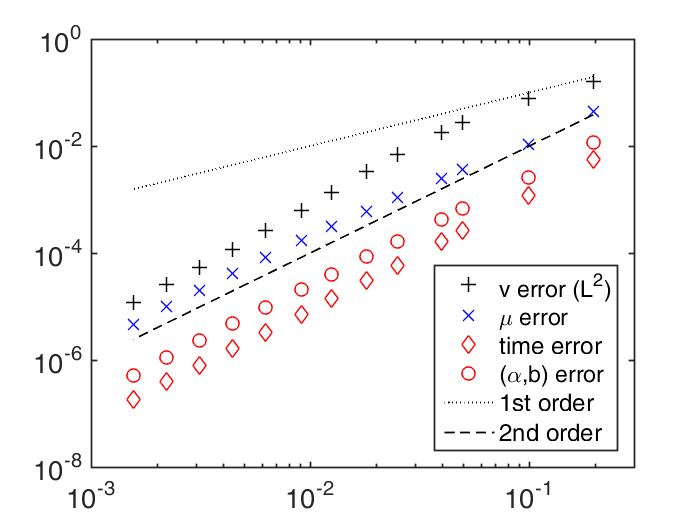}\label{fig:CTO}}
  \subfigure[$\nu=1$, fixed phase]{\includegraphics[height=5cm]%
    {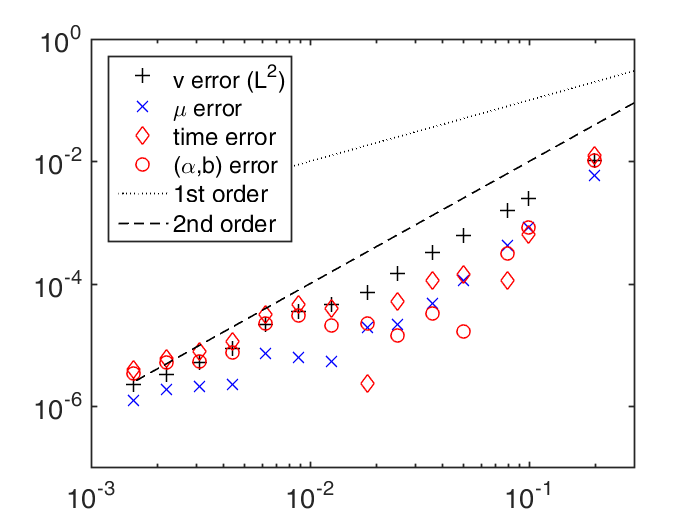}\label{fig:CLF}}
  \subfigure[$\nu=0.01$, fixed phase]{\includegraphics[height=5cm]%
    {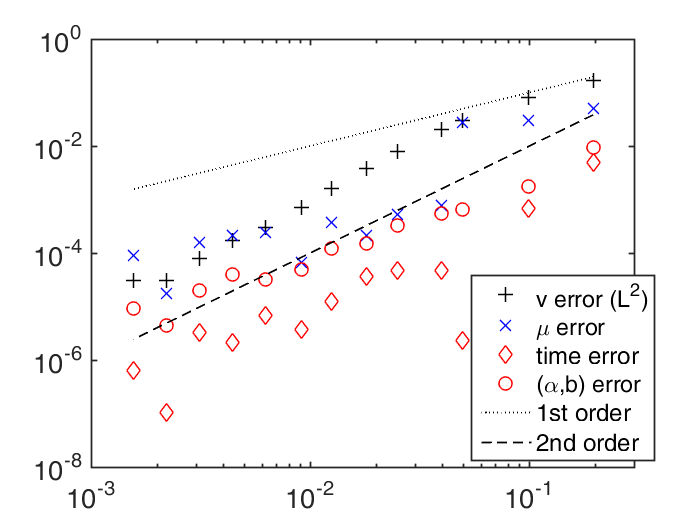}\label{fig:CTF}}
  \caption{Convergence plots for freezing 1d-Burgers equation.  With
    large viscosity left column)
    and very small viscosity (right column), and orthogonal phase
  condition (top row) and fixed phase condition (bottom row).}
  \label{fig:1dconvergence}
\end{figure} 
We calculate the error of all solution components separately.  Namely we
calculate the $L^2$-error of the differential
variable (``$v$ error'')
$\int_{R^-}^{R^+} |v_{\Delta\xi}(\xi,1)-v_\text{ref}(\xi,1)|^2\,d\xi$,
the error of the algebraic variables (``$\mu$ error'')
$|\mu_{\Delta\xi}(1)-\mu_\text{ref}(1)|_\infty$,
the error of the reconstructed time (``time error'')
$|t_{\Delta\xi}(1)-t_\text{ref}(1)|$,
and the error of the reconstructed transformation (``$\alpha,b$ error'')
$\max(|\alpha_{\Delta\xi}(1)-\alpha_\text{ref}(1)|,
|b_{\Delta\xi}(1)-b_\text{ref}(1)|)$.
Here a sub-index ${\Delta\xi}$ denotes the result of the numerical
approximation with step-size $\Delta\xi$ and a sub-index $\text{ref}$
refers to the reference solution.  The results are shown
in Fig.~\ref{fig:1dconvergence}, where we consider large and small
viscosities ($\nu=1$ and $\nu=0.01$) and in both cases orthogonal and
fixed phase conditions.  The results for other values of viscosity $\nu$
look very similar and are not shown.
In all experiments  $\Delta \xi$ is prescribed, and $\Delta \tau$
is related to $\Delta\xi$, so that \eqref{eq:deltatstep} holds, as
described above.  Moreover, we choose $R^+=-R^-=10$ for $\nu=1$ and 
$R^+=-R^-=5$ for $\nu=0.01$ so we may neglect the
influece of the boundary conditions.

On the horizontal axis in Fig.~\ref{fig:1dconvergence} we plot the
$\Delta \xi$-value and the vertical
axis is the error of the respective solution components.  In all
experiments one observes second order convergence for the
orthogonal phase condition, as was proved (for ODEs) in
Theorem~\ref{thm:2ndorder}.  One also observes second order convergence
for the differential variables in case of the fixed phase
condition, which was also shown in Theorem~\ref{thm:2ndorder}.

\textbf{Violation of CFL condition.}
Our final experiment for the 1d Burgers' equation concerns violation of
\eqref{eq:deltatstep}.
\begin{figure}[htb]
  \centering
  \subfigure[{$\nu=0$, $\lambda_\mathrm{CFL}=1.2$}]{\includegraphics[height=5cm]%
    {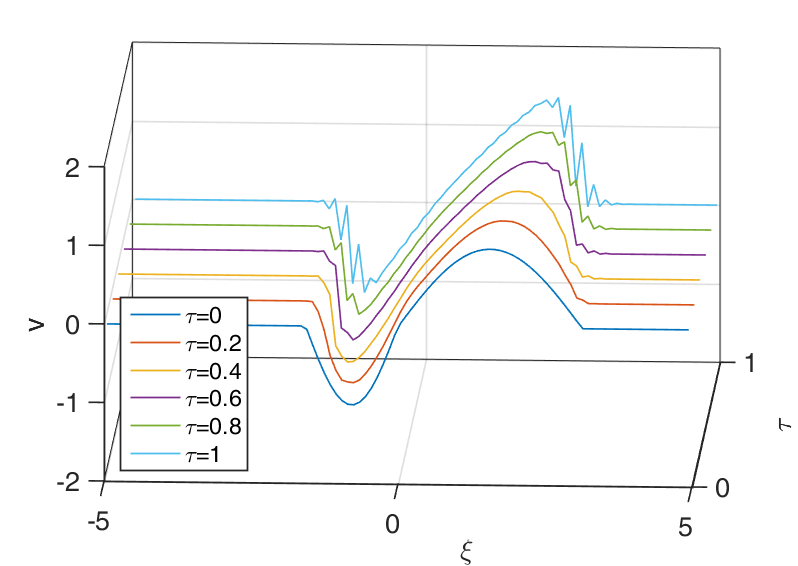}\label{fig:ECFL10a}}
    \subfigure[{$\nu=0.02$, $\lambda_\mathrm{CFL}=5$}]{\includegraphics[height=5cm]%
    {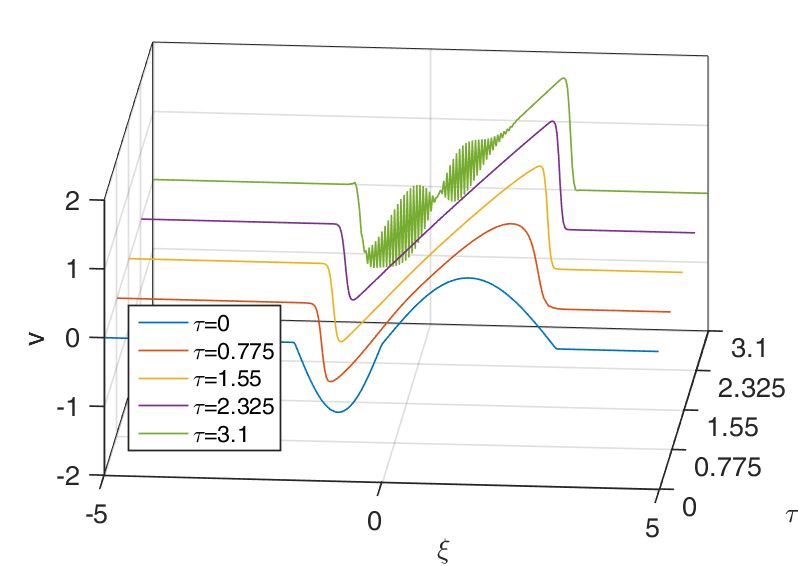}\label{fig:ECFL11a}}
    \caption{Appearance of oscillations, when the ratio
    $\tfrac{\Delta\tau}{\Delta\xi}$ is too large.}
\end{figure}
In Fig.~\ref{fig:ECFL10a} we consider the inviscid equation
and choose $\lambda_\mathrm{CFL}=1.2$ in \eqref{eq:deltatstep} (instead of
$\lambda_\mathrm{CFL}=\tfrac{1}{3}$).  The spatial stepsize is
$\Delta\xi=0.1$.
One nicely observes how oscillations develop.  Actually, the
calculation breaks down soon after the plotted solution at $\tau=1$.
For this choice of $\lambda_\mathrm{CFL}$ we observed the same behavior
also for all other step-sizes $\Delta\xi$
we have tried and we do not present these experiments here.

In Fig.~\ref{fig:ECFL11a} we show the result for $\nu=0.02$.  In
this case we have chosen $\lambda_\mathrm{CFL}=5$ in \eqref{eq:deltatstep}
and the spatial step-size was $\Delta\xi=10/350$.
Further experiments with different step-sizes $\Delta\xi$ 
seem to suggest that the numerical discretization of the
parabolic part is strongly smoothing and there is no linear barrier
of $(\Delta \tau)/(\Delta \xi)$ which
prevents oscillations but rather the ratio may increase as $\Delta \xi$
decreases.  Nevertheless, we expect that in the coupled
hyperbolic-parabolic case the linear relation dominates.

\subsection{2d-Experiments}
We also apply our method to the following two-dimensional Burgers' equations
\begin{equation}\label{eq:expBurgers}
  \partial_t u=\nu \Delta
  u-\tfrac{2}{3}(\partial_x+\partial_y)\bigl(|u|^\frac{3}{2}\bigr).
\end{equation}
That is, we numerically solve \eqref{eq:2dexpl} with the scheme
introduced in Sections~\ref{sec:3}--\ref{sec:4}.  First we choose a
fixed rectangular spatial domain $\Omega$, as explained in Section~\ref{sec:3.2}.
As in the one-dimensional case, we bound the maximal local wave speed by
the rough upper estimate
\begin{equation}\label{eq:wspeed2d}
  \mathbf{a}=\sup_{\xi\in\Omega}\sqrt{|v(\xi)|}
  +\max\bigl(|\mu_1|\max_{\xi\in\Omega}|\xi_1|+|\mu_2|,
  |\mu_1|\max_{\xi\in\Omega}|\xi_2|+|\mu_3|\bigr),
\end{equation}
and choose in each time step a step-size $\Delta\tau$ so that
\begin{equation}\label{eq:deltatstep2}
  \Delta \tau\le\frac{\min(\Delta \xi_1,\Delta \xi_2)}{
    \mathbf{a}}\cdot\lambda_\mathrm{CFL}.
\end{equation}
We found that $\lambda_\mathrm{CFL}=0.2$ is a suitable choice.  Note that in
\cite{KurganovTadmor:2000} there are bounds for $\lambda_\mathrm{CFL}$
given, which guarantee a maximum principle.

\begin{figure}[t!]
  \centering
  \subfigure[Initial data]{\includegraphics[height=4.2cm]%
    {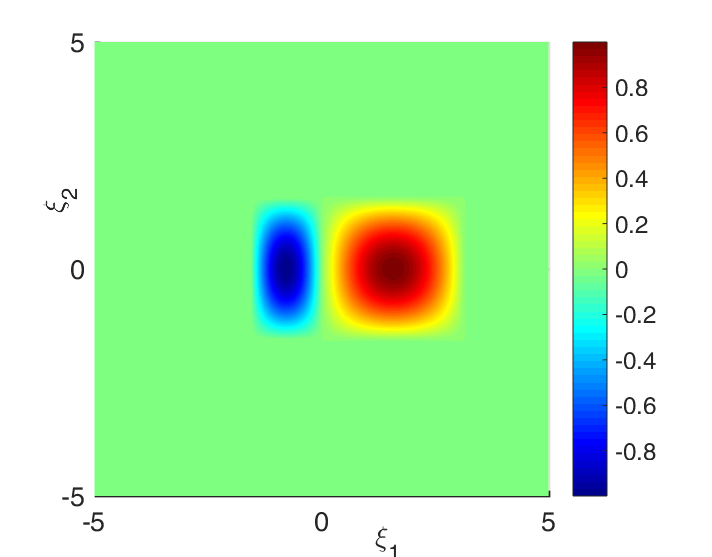}\label{fig:2dE1a}}%
  \subfigure[Final state at $\tau=6$]{\includegraphics[height=4.2cm]%
    {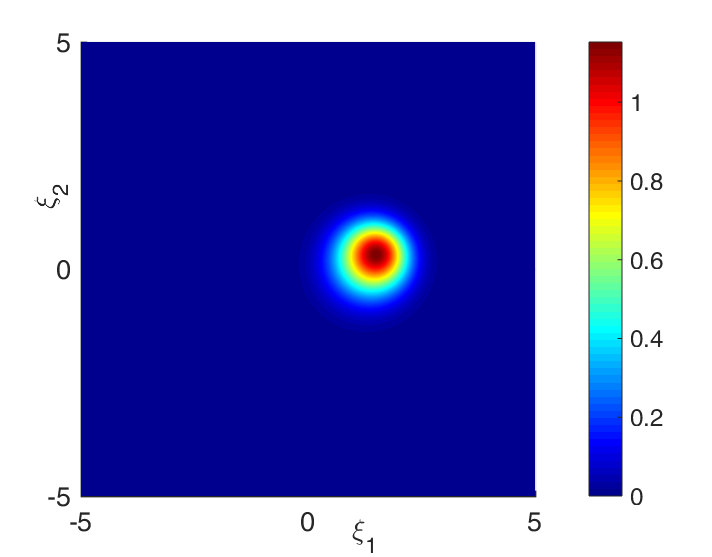}\label{fig:2dE1b}}%
  \subfigure[Relative error at $\tau=6$]{\includegraphics[height=4.2cm]%
    {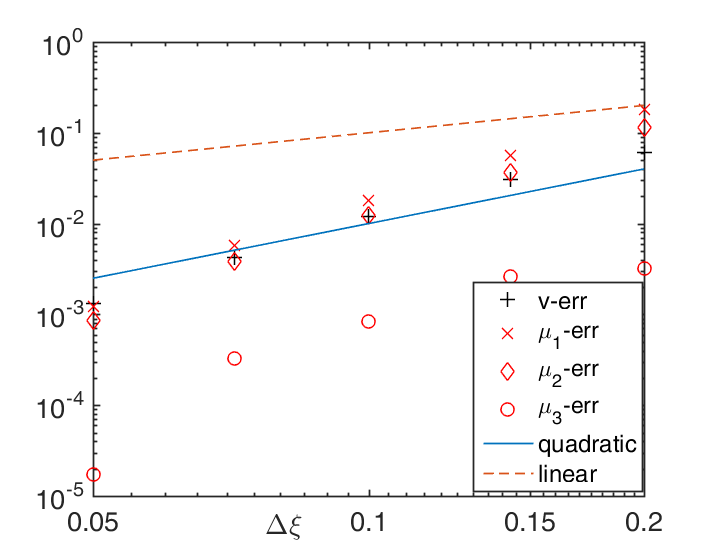}\label{fig:2dE1c}}
    \caption{The final state of the reference solution
      $v_\text{ref}(6)$ is shown in (b), it is obtained by the freezing
      method \eqref{eq:2dexpl} with $\nu=0.4$ and initial
      function shown in (a).  In (c) we show the
      relative errors of
      solutions obtained on coarser grids compared to the
      solution on the fine grid.}
  \label{fig:2dE1}
\end{figure}

Plots, showing the time evolution of the solution to the freezing
method, calculated with our scheme can be found in \cite{Rottmann:2016a}
and we do not repeat them here, as we focus on the discretization and
the error introduced by the discretization in the present article.  In
all our experiments we solve the Cauchy problem for \eqref{eq:burgersdd}
with initial data $u_0$, given by
\[u_0(x,y) = \begin{cases}
    \cos(y)\cdot\sin(2x),&-\tfrac{\pi}{2}<y<\tfrac{\pi}{2},-\tfrac{\pi}{2}<x<0,\\
    \cos(y)\cdot\sin(x),&-\tfrac{\pi}{2}<y<\tfrac{\pi}{2},0<x<\pi,\\
    0,&\text{otherwise},
\end{cases}\]
depicted in Fig.~\ref{fig:2dE1a}, by
the freezing method \eqref{eq:2dexpl}.  We choose the fixed domain
$\Omega=[-5,5]^2$ with no-flux boundary conditions.

In our first experiment we choose the viscosity $\nu=0.4$ and
solve the freezing system until $\tau=6$, which is a time when the
solution has settled.  The final state of this calculation obtained for
the fine grid with $\Delta\xi=\Delta\xi_1=\Delta\xi_2=\tfrac{1}{28}$ is shown in
Fig.~\ref{fig:2dE1b}.
We choose this solution as reference solution.  In
Fig.~\ref{fig:2dE1c} we plot the difference of the
solution components for different step-sizes.  More precisely we plot 
the relative $L^2$-difference 
$\tfrac{1}{\|v_\text{ref}\|}\|v_{\Delta\xi}-v_\text{ref}\|$
(``$v$-err''), and the relative differences of the algebraic variables
$\tfrac{1}{|\mu_{j,\text{ref}}|}|\mu_{j,\Delta\xi}-\mu_{j,\text{ref}}|$
(``$\mu_j$ error'') at the final time $\tau=6$.  The numerical findings
suggest that the scheme is in fact second order convergent (in
$\Delta\xi$).
\begin{figure}[htb]
  \centering
  \subfigure[Reference solution $v_\text{ref}$ at $\tau=300$ for cell
    size $0.05\times 0.05$.]{\includegraphics[height=4.2cm]%
    {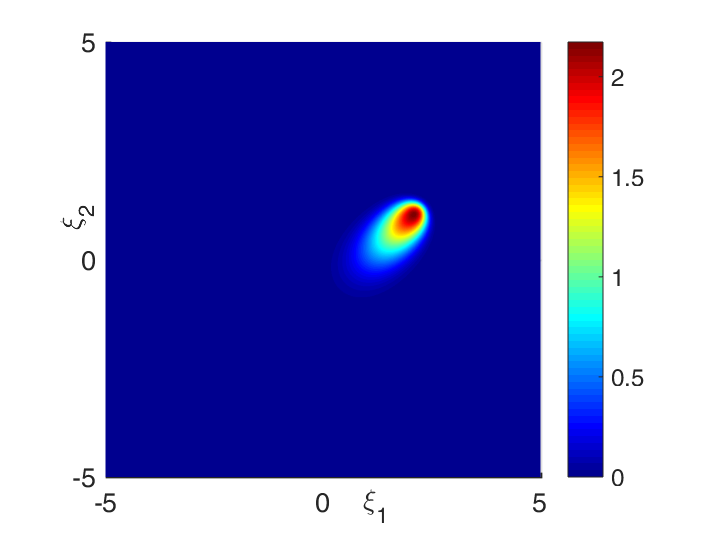}\label{fig:2dE2a}}%
    \subfigure[Difference $v_\text{ref}-v$.]{\includegraphics[height=4.2cm]%
    {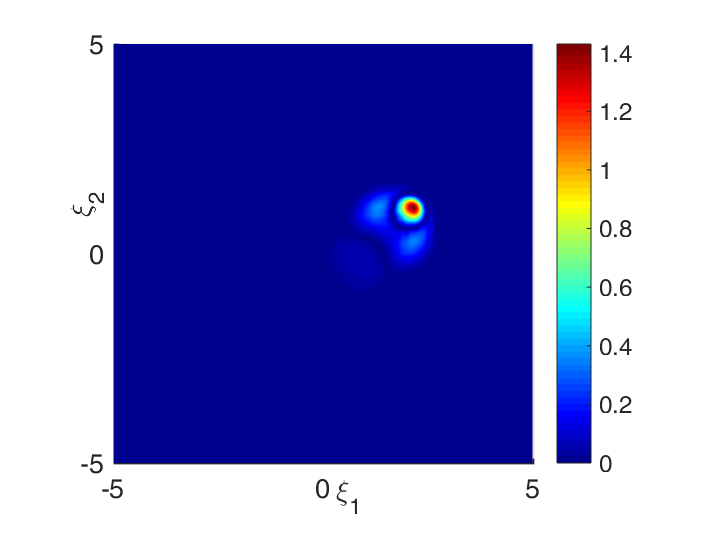}\label{fig:2dE2b}}%
  \subfigure[Difference $v_\text{ref}-v$.]{\includegraphics[height=4.2cm]%
    {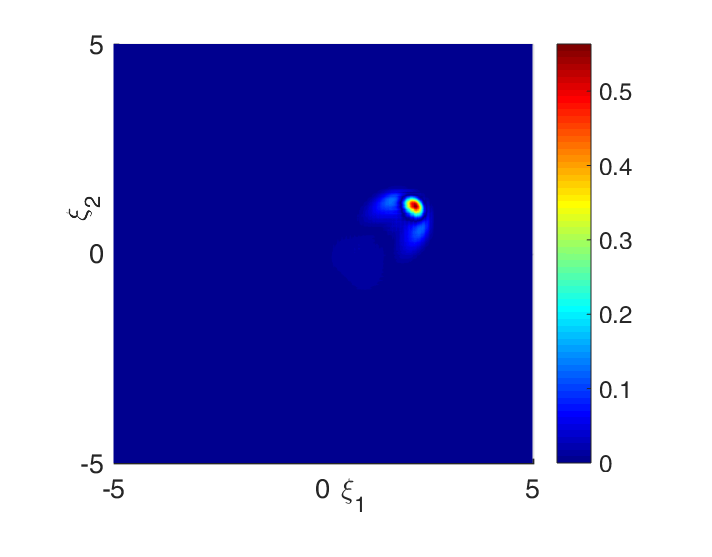}\label{fig:2dE2c}}
    \caption{Comparison of the solutions obtained with the freezing
      method at $\tau=300$ for $\nu=0.05$
      and different cell sizes, $0.2\times 0.2$ vs. $0.05\times 0.05$
    in (b), and $0.1\times 0.1$ vs. $0.05\times 0.05$ in (c).}
  \label{fig:2dE2}
\end{figure}
We repeat the experiment with the very small viscosity $\nu=0.05$.  The
solution to the 2d-Burgers' equation \eqref{eq:expBurgers} with this
viscosity and the same initial data as before shows a
meta-stable behavior similar to the one-dimensional case.  Therefore we
have to calculate for a very long time and we show the final
state of the calculation at time $\tau=300$ for
$\Delta\xi_1=\Delta\xi_2=0.05$ in Fig.~\ref{fig:2dE2a} and plot in
Fig.~\ref{fig:2dE2b},\ref{fig:2dE2c} the deviation of solutions
obtained for coarser grids with $\Delta\xi_1=\Delta\xi_2=0.2$ and
$\Delta\xi_1=\Delta\xi_2=0.1$ from the reference solution.  

\begin{figure}[htb]
  \centering
  \subfigure[Relative $L^2$-error depending on time.]{\includegraphics[height=4cm]%
    {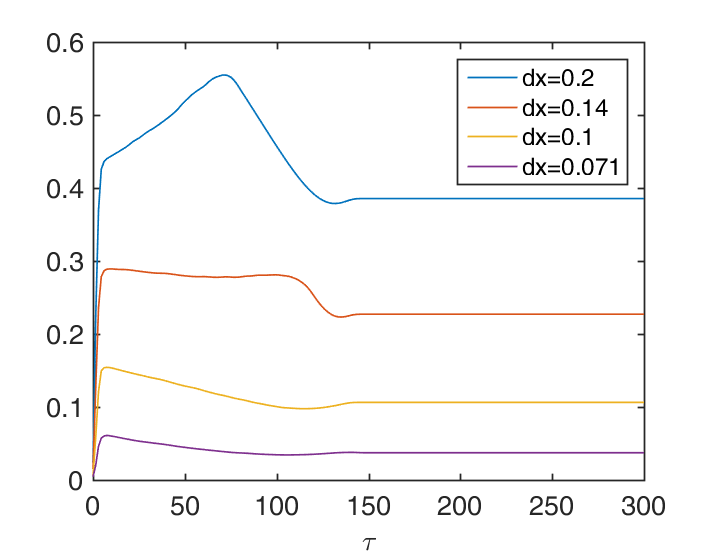}\label{fig:2dE3a}}%
  \subfigure[Zoom into initial behavior.]{\includegraphics[height=4cm]%
    {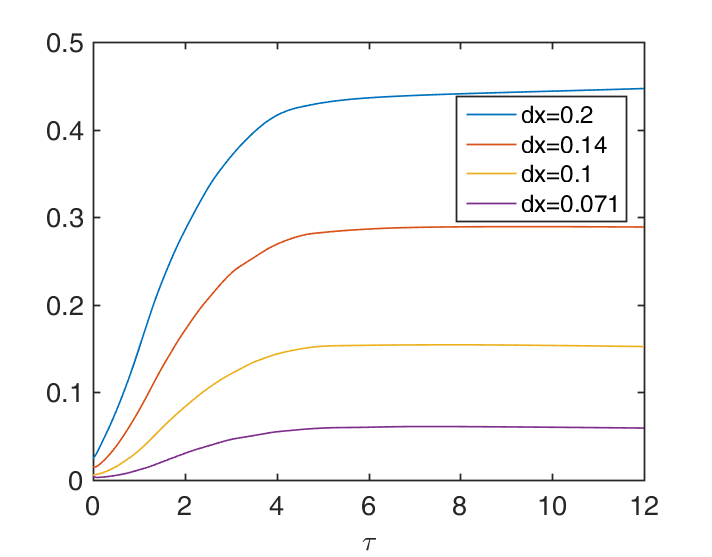}\label{fig:2dE3b}}
    \caption{The relative $L^2$-error
      $\tfrac{1}{\|v_\text{ref}(\tau)\|}\|v_\text{ref}(\tau)-v(\tau)\|$
      for different step-sizes $\Delta\xi_1=\Delta\xi_2=dx$ drawn as a
    function of time.}
  \label{fig:2dE3}
\end{figure}
\begin{figure}[htb]
  \centering
  \subfigure[Relative $\mu_1$-error depending on
  time.]{\includegraphics[height=4cm]%
    {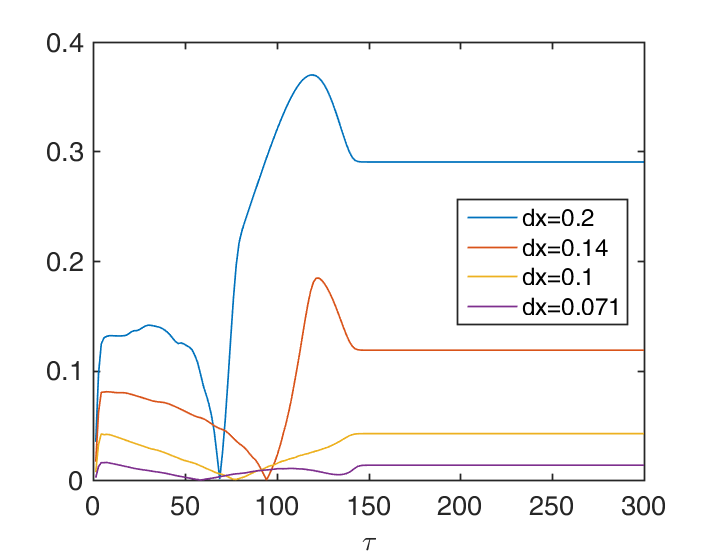}\label{fig:2dE4a}}%
  \subfigure[Zoom into initial behavior.]{\includegraphics[height=4cm]%
    {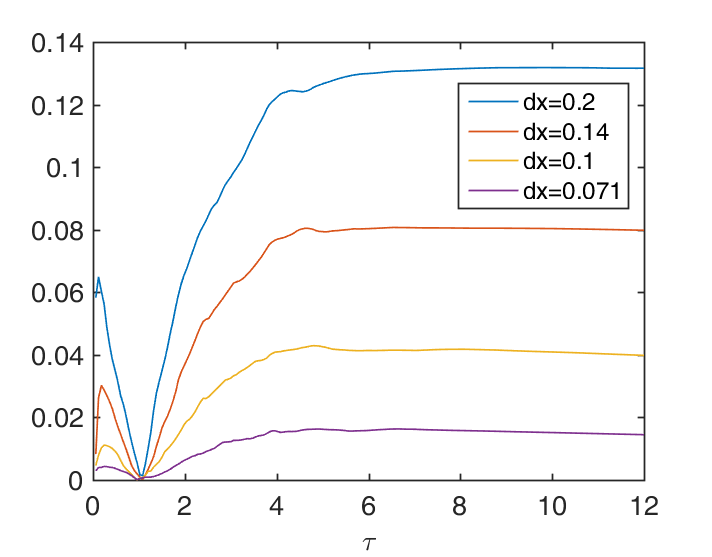}\label{fig:2dE4b}}
    \caption{The relative error
      $\tfrac{|\mu_{1,\text{ref}}(\tau)-\mu_1(\tau)|}{|\mu_{1,\text{ref}}(\tau)|}$
      of the algebraic variable corresponding
      to the speed of the scaling $\mu_1$
      for different step-sizes $\Delta\xi_1=\Delta\xi_2=dx$ drawn as a
    function of time.}
  \label{fig:2dE4}
\end{figure}
In Figs.~\ref{fig:2dE3}, \ref{fig:2dE4} we again consider the solution
to the freezing method for \eqref{eq:expBurgers} with $\nu=0.05$.  We
choose the solution with stepsize $\Delta\xi_1=\Delta\xi_2=0.05$ as
reference solution and plot for different
step-sizes $\Delta\xi_1=\Delta\xi_2=dx$, how the error of the
solution components depends on time.  For all solution components the
result is similar:  One observes that the error
initially rapidly increases, as it is expected since the error cumulates over
time.  But after a short time ($\tau\approx 6$) it settles and even
slowly decreases lateron ($\tau\approx 100$) until it finally reaches a
stationary value ($\tau\approx 150$).

\textbf{Violation of CFL condition.} We also observed that oscillations
appear, when we violate \eqref{eq:deltatstep2}.  We found that
oscillations for example appear for the same problem with $\nu=0.02$ and
$\lambda_{\mathrm{CFL}}=0.8$.  But we do not present the results here.



\textbf{Acknowledgement.}
We gratefully acknowledge financial support by the Deutsche
Forschungsgemeinschaft (DFG) through CRC 1173.\\
The author would like to thank Wolf-J\"urgen Beyn for many helpful
discussions and encouraging him to see this project through to the
finish.

\bibliographystyle{abbrv}

\begin{appendix}
\section{Discretization of the hyperbolic part}
\label{sec:KT}
For the discretization of the hyperbolic part
in~\eqref{eq:CoCauchyA1}, we
briefly review and adapt a scheme introduced by Kurganov and Tadmor in
\cite{KurganovTadmor:2000} (KT-scheme).  In
particular, we need an adaptation for non-homogeneous flux functions, i.e.\
the case of hyperbolic conservation laws of the form
\[u_\tau+\frac{d}{d\xi}f(\xi,u)=0.\]
In this appendix $u$ may be $\R^m$ valued, but in the actual application in
Section~\ref{sec:3} it is always a scalar.
For the KT-scheme, a uniform grid $\xi_j=\xi_0+j\Dxi$,
$j\in\Z$ in $\R$ is chosen and one assumes that at a time instance $\tau^n$
an approximation of the average mass of the function $u$ in the spatial cells
$(\xi_{j-\halb},\xi_{j+\halb})$, $\xi_{j\pm\halb}=\xi_j\pm\tfrac{\Delta\xi}{2}$
$j\in\Z$, is given as $u_j^n\in\R^m$.  I.e.\ $u_j^n\approx
\frac{1}{\Dxi}\int_{\xi_{j-\halb}}^{\xi_{j+\halb}} u(\xi,t^n)\,d\xi$.  From
this grid function one obtains the piecewise linear reconstruction at the time
instance $\tau^n$ as
\[\wt{u}(\xi,t^n)=
\sum_j\left( u_j^n\mathbbm{1}_{[\xi_{j-\halb},\xi_{j+\halb})}(\xi)+
  (u_\xi)_j^n(\xi-\xi_j)\mathbbm{1}_{[\xi_{j-\halb},\xi_{j+\halb})}(\xi)\right),\]
where $(u_\xi)_j^n$ is a suitable choice of the slope of this piecewise
linear function in the cell $(\xi_{j-\halb},\xi_{j+\halb})$.  These ``suitable''
slopes are given by the ``minmod''-reconstruction, i.e. for a fixed
$\theta\in[1,2]$ the slopes are given as
\[(u_\xi)_j^n = \mm\left(\theta \frac{u_j-u_{j-1}}{\Dxi},
  \frac{u_{j+1}-u_{j-1}}{2\Dxi},
  \theta\frac{u_{j+1}-u_j}{\Dxi}\right),\]
where
\begin{equation}\label{eq:app01}
  \mm(a_1,\ldots,a_k)=\max\Bigl(\min(a_1,0),\ldots,\min(a_k,0)\Bigr)+
  \min\Bigl(\max(a_1,0),\ldots,\max(a_k,0)\Bigr). 
\end{equation}
In the vector-valued case \eqref{eq:app01} is understood in the component wise
sense.
The idea of the KT-scheme is to use an artificial finer grid, which depends on the
time step size $\Delta\tau$ and separates the smooth and non-smooth regions of the
solution.  In the original paper \cite{KurganovTadmor:2000}, local
(maximal) wave speeds are calculated and a spatial grid based on these is
chosen.  Since we ultimately intend to apply the scheme to the
PDAE-system~\eqref{eq:CoCauchy}, for which the wave speeds nonlinearly depend on
the solution, we choose a sufficiently large
$a>0$, which bounds the spectral radius of $\tfrac{\partial}{\partial
  u} f(\xi ,u)$.  Note that this is obviously not possible for
$f(\xi ,u)=\xi u$, $\xi \in\R$, but because we will restrict to compact domains
in the end, we assume that the spectral radius is uniformly bounded in
$\xi $ for bounded $u$.

The new grid is then given by
\[\ldots<\xi_{j-1}<\xi_{j-\halb,l}=\xi_{j-\halb}-a\Dtau<\xi_{j-\halb,r}=\xi_{j-\halb}+a\Dtau<\xi_j<\ldots\]
for $\Dtau$ sufficiently small.  We remark that for $\Dtau\searrow 0$, the new
grid points $\xi_{j-\halb,l}$ and $\xi_{j-\halb,r}$ both converge to
$\xi_{j-\halb}$, the first from the left and the second from the right.
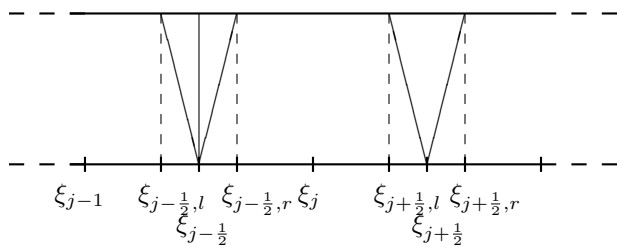
\begin{figure}[htb]
  \centering
  \setlength{\unitlength}{4144sp}%
\begingroup\makeatletter\ifx\SetFigFont\undefined%
\gdef\SetFigFont#1#2#3#4#5{%
  \reset@font\fontsize{#1}{#2pt}%
  \fontfamily{#3}\fontseries{#4}\fontshape{#5}%
  \selectfont}%
\fi\endgroup%
\begin{picture}(3644,1382)(-21,-521)
\thicklines
{\color[rgb]{0,0,0}\multiput(  1,839)(180.00000,0.00000){3}{\line( 1, 0){ 90.000}}
}%
{\color[rgb]{0,0,0}\multiput(  1,-61)(180.00000,0.00000){3}{\line( 1, 0){ 90.000}}
}%
{\color[rgb]{0,0,0}\multiput(3151,839)(180.00000,0.00000){3}{\line( 1, 0){ 90.000}}
}%
{\color[rgb]{0,0,0}\multiput(3151,-61)(180.00000,0.00000){3}{\line( 1, 0){ 90.000}}
}%
{\color[rgb]{0,0,0}\put(361,-61){\line( 1, 0){2790}}
}%
{\color[rgb]{0,0,0}\put(3151,839){\line(-1, 0){2790}}
}%
{\color[rgb]{0,0,0}\put(451,-16){\line( 0,-1){ 90}}
}%
{\color[rgb]{0,0,0}\put(901,-16){\line( 0,-1){ 90}}
}%
{\color[rgb]{0,0,0}\put(1126,-16){\line( 0,-1){ 90}}
}%
{\color[rgb]{0,0,0}\put(1351,-16){\line( 0,-1){ 90}}
}%
{\color[rgb]{0,0,0}\put(1801,-16){\line( 0,-1){ 90}}
}%
{\color[rgb]{0,0,0}\put(2251,-16){\line( 0,-1){ 90}}
}%
{\color[rgb]{0,0,0}\put(2476,-16){\line( 0,-1){ 90}}
}%
{\color[rgb]{0,0,0}\put(2701,-16){\line( 0,-1){ 90}}
}%
{\color[rgb]{0,0,0}\put(3151,-16){\line( 0,-1){ 90}}
}%
\put(271,-286){\makebox(0,0)[lb]{\smash{{\SetFigFont{10}{12.0}{\familydefault}{\mddefault}{\updefault}{\color[rgb]{0,0,0}$\xi_{j-1}$}%
}}}}
\put(766,-286){\makebox(0,0)[lb]{\smash{{\SetFigFont{10}{12.0}{\familydefault}{\mddefault}{\updefault}{\color[rgb]{0,0,0}$\xi_{j-\halb,l}$}%
}}}}
\put(991,-466){\makebox(0,0)[lb]{\smash{{\SetFigFont{10}{12.0}{\familydefault}{\mddefault}{\updefault}{\color[rgb]{0,0,0}$\xi_{j-\halb}$}%
}}}}
\put(1261,-286){\makebox(0,0)[lb]{\smash{{\SetFigFont{10}{12.0}{\familydefault}{\mddefault}{\updefault}{\color[rgb]{0,0,0}$\xi_{j-\halb,r}$}%
}}}}
\put(1711,-286){\makebox(0,0)[lb]{\smash{{\SetFigFont{10}{12.0}{\familydefault}{\mddefault}{\updefault}{\color[rgb]{0,0,0}$\xi_j$}%
}}}}
\put(2386,-466){\makebox(0,0)[lb]{\smash{{\SetFigFont{10}{12.0}{\familydefault}{\mddefault}{\updefault}{\color[rgb]{0,0,0}$\xi_{j+\halb}$}%
}}}}
\put(2161,-286){\makebox(0,0)[lb]{\smash{{\SetFigFont{10}{12.0}{\familydefault}{\mddefault}{\updefault}{\color[rgb]{0,0,0}$\xi_{j+\halb,l}$}%
}}}}
\put(2611,-286){\makebox(0,0)[lb]{\smash{{\SetFigFont{10}{12.0}{\familydefault}{\mddefault}{\updefault}{\color[rgb]{0,0,0}$\xi_{j+\halb,r}$}%
}}}}
\thinlines
{\color[rgb]{0,0,0}\multiput(901,-61)(0.00000,120.00000){8}{\line( 0, 1){ 60.000}}
}%
{\color[rgb]{0,0,0}\multiput(1351,-61)(0.00000,120.00000){8}{\line( 0, 1){ 60.000}}
}%
{\color[rgb]{0,0,0}\multiput(2251,-61)(0.00000,120.00000){8}{\line( 0, 1){ 60.000}}
}%
{\color[rgb]{0,0,0}\multiput(2701,-61)(0.00000,120.00000){8}{\line( 0, 1){ 60.000}}
}%
{\color[rgb]{0,0,0}\put(901,839){\line( 1,-4){225}}
}%
{\color[rgb]{0,0,0}\multiput(1126,-61)(-1.31223,7.87340){115}{\makebox(1.5875,11.1125){\tiny}}
}%
{\color[rgb]{0,0,0}\put(1126,-61){\line( 0, 1){900}}
}%
{\color[rgb]{0,0,0}\multiput(1126,-61)(1.31308,7.87850){114}{\makebox(1.5875,11.1125){\tiny}}
}%
{\color[rgb]{0,0,0}\put(1126,-61){\line( 1, 4){225}}
}%
{\color[rgb]{0,0,0}\put(2476,-61){\line(-1, 4){225}}
}%
{\color[rgb]{0,0,0}\multiput(2476,-61)(-1.31308,7.87850){114}{\makebox(1.5875,11.1125){\tiny}}
}%
{\color[rgb]{0,0,0}\multiput(2476,-61)(1.31308,7.87850){114}{\makebox(1.5875,11.1125){\tiny}}
}%
{\color[rgb]{0,0,0}\put(2476,-61){\line( 1, 4){225}}
}%
\end{picture}%
  \caption{Sketch of the grid and smooth and non-smooth parts of the solution.}
  \label{fig:KT}
\end{figure}
The principal idea now is, to use the integral form of the conservation law in
the smooth and non-smooth regions for the calculation of the mass at the new time
instance $\tau^{n+1}$.  We begin with the non-smooth part:
\begin{multline*}
  \frac{1}{\Dxihalb} \int_{\xi_{j-\halb,l}}^{\xi_{j-\halb,r}}
  u(\xi,\tau^{n+1})\,d\xi = \frac{1}{2a\Dtau}\Bigl(
  \int_{\xi_{j-\halb,l}}^{\xi_{j-\halb,r}} \wt{u}(\xi,\tau^n)\,d\xi\\
  -\int_{\tau^n}^{\tau^{n+1}}
  f\bigl(\xi_{j-\halb,r},u(\xi_{j-\halb,r},\tau)\bigr)\,d\tau
  +\int_{\tau^n}^{\tau^{n+1}}
  f\bigl(\xi_{j-\halb,l},u(\xi_{j-\halb,l},\tau)\bigr)\,d\tau\Bigr),
\end{multline*}
where $\Dxihalb=2a\Dtau=\xi_{j-\halb,r}-\xi_{j-\halb,l}$.
The time integrals are in the smooth regions of the solution (the
dashed vertical lines in Figure~\ref{fig:KT}).  We approximate them by the
midpoint rule, so that we need the value of the solution $u$ at
$\xi_{j-\halb,l}$, $\xi_{j-\halb,r}$ and time instance $\tau^n+\tfrac{\Dtau}{2}$.  In a smooth
region this value of $u$ is approximately
\begin{align*}
  u(\xi_{j-\halb,l},\tau^{n+\halb})&\approx
  u_{j-\halb,l}^n-\tfrac{\Dtau}{2}\left(
    \tfrac{\partial}{\partial \xi}f(\xi_{j-\halb,l},u_{j-\halb,l}^n)+
    \tfrac{\partial}{\partial
      u}f(\xi_{j-\halb,l},u_{j-\halb,l}^n)(u_\xi)_j^n\right)\\
  &:=u_{j-\halb,l}^{n+\halb}
\end{align*}
and the same for $l$ replaced by $r$.  For the average value of the solution $u$
at the new time instance $\tau^{n+1}$ in the spatial interval
$(\xi_{j-\halb,l},\xi_{j+\halb,r})$ this yields the approximation
\begin{equation}\label{eq:wjhalb}
  \begin{aligned}
    w_{j-\halb}^{n+1}&:=\frac{u_{j-1}^n+u_j^n}{2}+ \frac{\Dxi-a\Dt}{4}
    \bigl((u_\xi)_{j-1}^n-(u_\xi)_j^n\bigr)
    -\frac{1}{2a}
    \bigl(f(\xi_{j-\halb,r},u_{j-\halb,r}^{n+\halb})-f(\xi_{j-\halb,l},u_{j-\halb,l}^{n+\halb})\bigr)\\
    &\approx \frac{1}{\Dxihalb} \int_{\xi_{j-\halb,l}}^{\xi_{j-\halb,r}}
    u(\xi,t^{n+1})\,d\xi.
  \end{aligned}
\end{equation}
Similarly, the average of
the solution at the new time instance $\tau^{n+1}$ in the spatial interval
$(\xi_{j-\halb,r},\xi_{j+\halb,l})$ is approximated by
\begin{equation}\label{eq:wj}
  \begin{aligned}
    w_j^{n+1}&:=u_j^n-\frac{\Dtau}{\Dxi-2a\Dtau}\bigl(f(\xi_{j+\halb,l},u_{j+\halb,l}^{n+\halb})-f(\xi_{j-\halb,r},u_{j-\halb,r}^{n+\halb})\bigr)\\
    &\approx
    \frac{1}{\Dxi-\Dxihalb}\int_{\xi_{j-\halb,r}}^{\xi_{j+\halb,l}}u(\xi,t^{n+1})\,d\xi.
  \end{aligned}
\end{equation}
To come back from the new grid back to the original one, a piecewise
linear reconstruction $\wt{w}$ of the grid function $w^{n+1}$ is
calculated and then $u_j^{n+1}$ is obtained by integration over the
cells $(\xi_{j-\halb},\xi_{j+\halb})$.  Finally, one obtains the fully
discrete scheme:
\begin{equation}\label{eq:fullyDiscreteKT}
  \begin{aligned}
    u_j^{n+1}&=\frac{1}{\Dxi} \Bigl[ (\xi_{j-\halb,r}-\xi_{j-\halb})
    w_{j-\halb}^{n+1}+
    (w_\xi)_{j-\halb}^{n+1} \frac{(\xi_{j-\halb,r}-\xi_{j-\halb})^2}{2}\\
    &\qquad\quad +(\xi_{j+\halb,l}-\xi_{j-\halb,r}) w_j^{n+1}\\
    &\qquad\quad +(\xi_{j+\halb}-\xi_{j+\halb,l})w_{j+\halb}^{n+1}
    - (w_\xi)_{j+\halb}^{n+1}
    \frac{(\xi_{j+\halb}-\xi_{j+\halb,l})^2}{2}\Bigr],
  \end{aligned}
\end{equation}
where $w_{j\pm\halb}^{n+1}$ and $w_j^{n+1}$ are given by
\eqref{eq:wjhalb} and \eqref{eq:wj} and
\[(w_\xi)_{j-\halb}^{n+1}=\frac{2}{\Dxi} \mm\bigl(
w_j^{n+1}-w_{j-\halb}^{n+1},w_{j-\halb}^{n+1}-w_{j-1}^{n+1}\bigr).\]
The fully discrete scheme~\eqref{eq:fullyDiscreteKT} admits a
semi-discrete version, i.e. a method of lines system, which we will use in
our discretization of \eqref{eq:CoCauchy}.  To see this, let
$\lambda=\Dtau/\Dxi$, in which we consider $\Dxi$ as a fixed quantity and are interested
in the limit $\Dtau\to 0$ so that
$\mathcal{O}(\lambda^2)=\mathcal{O}(\Dtau^2)$.
We note, that the summands in~\eqref{eq:fullyDiscreteKT} can be written in the following
form, where we use
$\xi_{j-\halb}-\xi_{j-\halb,l}=\xi_{j-\halb,r}-\xi_{j-\halb}=a\Dtau$:
\begin{align}\label{eq:fullyDiscI}
  (\xi_{j-\halb,r}-\xi_{j-\halb}) w_{j-\halb}^{n+1}&=\frac{a\Dtau}{2}
  \bigl(u_{j-1}^n+\frac{\Dxi}{2}(u_\xi)_{j-1}^n+u_j^n-\frac{\Dxi}{2}
  (u_\xi)_j^n\bigr)\\\notag
  &\quad-\frac{\Dtau}{2}\bigl(f(\xi_{j-\halb,r},u_{j-\halb,r}^{n+\halb})-f(\xi_{j-\halb,l},u_{j-\halb,l}^{n+\halb})\bigr)
  +\mathcal{O}(\lambda^2),\\
  \label{eq:fullyDiscII}
  (w_\xi)_{j-\halb}^{n+1} \frac{(\xi_{j-\halb,r}-\xi_{j-\halb})^2}{2}&=\mathcal{O}(\lambda^2)\\
  \label{eq:fullyDiscIII}
  (\xi_{j+\halb,l}-\xi_{j-\halb,r}) w_j^{n+1}&=\Dxi u_j^n - 2a\Dtau u_j^n-\Dtau
  \bigl( f(\xi_{j+\halb,l},u_{j+\halb,l}^{n+\halb})-f(\xi_{j-\halb,r},u_{j-\halb,r}^{n+\halb})\bigr),\\
  \label{eq:fullyDiscIV}
  (\xi_{j+\halb}-\xi_{j+\halb,l})w_{j+\halb}^{n+1}&=\frac{a\Dtau}{2}
    \bigl(u_{j}^n+\frac{\Dxi}{2}(u_\xi)_{j}^n+u_{j+1}^n-\frac{\Dxi}{2}
    (u_\xi)_{j+1}^n\bigr)\\\notag
    &\quad-\frac{\Dtau}{2}\bigl(f(\xi_{j+\halb,r},u_{j+\halb,r}^{n+\halb})-f(\xi_{j+\halb,l},u_{j+\halb,l}^{n+\halb})\bigr)
    +\mathcal{O}(\lambda^2),\\
  \label{eq:fullyDiscV}
  (w_\xi)_{j+\halb}^{n+1}
  \frac{(\xi_{j+\halb}-\xi_{j+\halb,l})^2}{2}&=\mathcal{O}(\lambda^2).
\end{align}
The semi-discrete version is then obtained by subtracting $u_j^n$ from
both sides of~\eqref{eq:fullyDiscreteKT}, dividing by $\Dtau$ and
considering the limit $\Dtau\to 0$.
Using~\eqref{eq:fullyDiscI}--\eqref{eq:fullyDiscV} and the
convergences $\xi_{j\pm\halb,l/r}\to \xi_{j\pm\halb}$,
$u_{j-\halb,l}^{n+\halb}\to
u_{j-\halb}^{n,-}=u_{j-1}^n+\frac{\Dxi}{2}(u_\xi)_{j-1}^n$, and
$u_{j-\halb,r}^{n+\halb}\to
u_{j-\halb}^{n,+}=u_{j}^n-\frac{\Dxi}{2}(u_\xi)_{j}^n$ as $\Dtau\to 0$, the
continuity properties of $f$ imply that the limit of the difference
quotient exists.  The final result is the following semi-discrete scheme (method of lines
system):
\begin{equation}\label{eq:semidiscKT}
  \begin{aligned}
    \frac{d}{d\tau} u_j&=
    \frac{1}{2\Dxi}\Bigl(\bigl(f(\xi_{j-\halb},u_{j-\halb}^-)+f(\xi_{j-\halb},u_{j-\halb}^+)\bigr)
    -\bigl(f(\xi_{j+\halb},u_{j+\halb}^-)+f(\xi_{j+\halb},u_{j+\halb}^+)\bigr)\Bigr)\\
    &\qquad +\frac{a}{2\Dxi}
    \bigl(u_{j-\halb}^--u_{j-\halb}^+-u_{j+\halb}^-+u_{j+\halb}^+\bigr),\quad
    j\in\Z,
  \end{aligned}
\end{equation}
where 
\begin{equation}\label{eq:minmodlim1d}
\begin{aligned}
  u_{j-\halb}^-&=u_{j-1}+\halb\mm\bigl(\theta
  (u_{j-1}-u_{j-2}),\frac{u_j-u_{j-2}}{2},\theta
  (u_j-u_{j-1})\bigr),\\
  u_{j-\halb}^+&=u_{j}-\halb\mm\bigl(\theta
  (u_{j}-u_{j-1}),\frac{u_{j+1}-u_{j-1}}{2},\theta
  (u_{j+1}-u_{j})\bigr).
\end{aligned}
\end{equation}
Before we continue, we state some simple remarks:
\begin{remark}\label{rem:KT}
  \begin{enumerate}
  \item In the case $f(\xi,u)=f(u)$,~\eqref{eq:semidiscKT} reduces to the
    semi-discrete version of the KT-scheme from
    \cite{KurganovTadmor:2000}.
  \item In the case of a
    piecewise constant reconstruction of $u$, the last summand in
    \eqref{eq:semidiscKT} reduces to
    $\tfrac{a}{2\Dxi}(u_{j-1}-2u_j+u_{j+1})$, and can be interpreted as
    artificial viscosity, that vanishes as $\Dxi\to 0$.
  \item The flux function $f$ enters linearly into the semi-discrete
    scheme, this is an important observation, because it allows us to
    derive a method of lines system for the PDE-part of the
    PDAE~\eqref{eq:CoCauchy}, without knowing the algebraic variables $\mu$ in
    advance.  To the resulting system we may apply
    a so called half-explicit scheme for the PDAE (see \cite{HairerLubichRoche:1989})
  \item The method of lines system \eqref{eq:semidiscKT} can be
    written in conservative form as
    \[ u_j'=-\frac{H_{j+\halb}-H_{j-\halb}}{\Dxi},\]
    where
    \[H_{j+\halb}=\frac{f(\xi_{j+\halb},u_{j+\halb}^+)+f(\xi_{j+\halb},u_{j+\halb}^-)}{2}
    - a\bigl(u_{j+\halb}^+-u_{j+\halb}^-\bigr).\]
  \end{enumerate}
\end{remark}
\textbf{Multi-dimensional version}\\
There is no difficulty in performing the same calculation in the
multi-dimensional setting.  We only present the result for the
2-dimensional case:  For this we consider a system of hyperbolic
conservation laws of the form
\[u_\tau+\frac{d}{d\xi}f(\xi,\eta,u)+
\frac{d}{d\eta }g(\xi,\eta,u)=0.\]
Choose a uniform spatial grid
$(\xi_j,\eta_k)=(\xi_0+j\Dxi,\eta_0+k\Deta)$ and assume that a grid function
$(u_{j,k})$ is given on this grid, where $u_{j,k}(\tau)$ denotes the
average mass in the cell $(\xi_{j-\halb},\xi_{j+\halb})\times
(\eta_{k-\halb},\eta_{k+\halb})$, $\xi_{j\pm\halb}=\xi_j\pm\tfrac{\Dxi}{2}$ and
similarly $\eta_{k\pm\halb}=\eta_k\pm\tfrac{\Deta}{2}$ at time $\tau$.
Adapting the above derivation, we obtain the following 2-dimensional
version of the semi-discrete KT-scheme:
\begin{equation}\label{eq:semidiscKT2d}
  \begin{aligned}
  u_{j,k}'&=\frac{1}{2\Dxi} \Bigl[ 
  \bigl(f(\xi_{j-\halb},\eta_k,u_{j-\halb,k}^-)+
  f(\xi_{j-\halb},\eta_k,u_{j-\halb,k}^+)\bigr)
  -\bigl(f(\xi_{j+\halb},\eta_k,u_{j+\halb,k}^-)+
  f(\xi_{j+\halb},\eta_k,u_{j+\halb,k}^+)\bigr)\Bigr]\\
  &\quad+\frac{1}{2\Deta} \Bigl[
  \bigl(g(\xi_j,\eta_{k-\halb},u_{j,k-\halb}^-)
  +g(\xi_j,\eta_{k-\halb},u_{j,k-\halb}^+)\bigr)
  -\bigl( g(\xi_j,\eta_{k+\halb},u_{j,k+\halb}^-)
  +g(\xi_j,\eta_{k+\halb},u_{j,k+\halb}^+)\bigr)\Bigr]\\
  &\quad+\frac{a}{2\Dxi}
  \Bigl(u_{j-\halb,k}^--u_{j-\halb,k}^+-u_{j+\halb,k}^-+
  u_{j+\halb,k}^+\Bigr)\\
  &\quad+\frac{a}{2\Deta}
  \Bigl(u_{j,k-\halb}^--u_{j,k-\halb}^+-
  u_{j,k+\halb}^-+u_{j,k+\halb}^+\Bigr).
\end{aligned}
\end{equation}
In~\eqref{eq:semidiscKT2d} the number $a$ is again assumed to be an
upper bound for the maximum of the spectral radii of
$\tfrac{\partial}{\partial u} f(\xi,\eta,u)$ and $\tfrac{\partial}{\partial
  u} g(\xi,\eta,u)$ for all relevant $\xi,\eta,u$.  We assume that such an upper
bound exists.

The one-sided limits appearing in~\eqref{eq:semidiscKT2d} are given by
\begin{equation}\label{eq:minmodlim2d}
\begin{aligned}
u_{j-\halb,k}^-&=u_{j-1,k}+\tfrac{1}{2}\mm\bigl(\theta(u_{j-1,k}-u_{j-2,k}),\frac{u_{j,k}-u_{j-2,k}}{2},\theta(u_{j,k}-u_{j-1,k})\bigr),\\
u_{j-\halb,k}^+&=u_{j,k}-\tfrac{1}{2}\mm\bigl(\theta(u_{j,k}-u_{j-1,k}),\frac{u_{j+1,k}-u_{j-1,k}}{2},\theta(u_{j+1,k}-u_{j,k})\bigr),\\
u_{j,k-\halb}^-&=u_{j,k-1}+\tfrac{1}{2}\mm\bigl(\theta(u_{j,k-1}-u_{j,k-2}),\frac{u_{j,k}-u_{j,k-2}}{2},\theta(u_{j,k}-u_{j,k-1})\bigr),\\
u_{j,k-\halb}^+&=u_{j,k}-\tfrac{1}{2}\mm\bigl(\theta(u_{j,k}-u_{j,k-1}),\frac{u_{j,k+1}-u_{j,k-1}}{2},\theta(u_{j,k+1}-u_{j,k})\bigr).
\end{aligned}
\end{equation}


\end{appendix}

\end{document}